\documentclass{amsart}

\usepackage{graphicx}

\newtheorem{theorem}{Theorem}[section]
\newtheorem{lemma}[theorem]{Lemma}
\newtheorem{proposition}[theorem]{Proposition}
\newtheorem{corollary}[theorem]{Corollary}

\newtheorem{_conjecture}[theorem]{Conjecture}

\newtheorem{_problem}[theorem]{Problem}

\newtheorem{_algorithm}[theorem]{Algorithm}

\newtheorem{_subroutime}[theorem]{Subroutine}

\newtheorem{_claim}[theorem]{Claim}
\newenvironment{claim}{\begin{_claim}\rm}{\end{_claim}}

\newtheorem{_subclaim}[theorem]{Sub-claim}

\newtheorem{_definition}[theorem]{Definition}
\newenvironment{definition}{\begin{_definition}\rm}{\end{_definition}}

\newtheorem{_remark}[theorem]{\it Remark}
\newenvironment{remark}{\begin{_remark}\rm}{\end{_remark}}

\newtheorem{_example}[theorem]{Example}
\newenvironment{example}{\begin{_example}\rm}{\end{_example}}

\newtheorem{_question}[theorem]{Question}
\newenvironment{question}{\begin{_question}\rm}{\end{_question}}

\numberwithin{equation}{section}
\numberwithin{table}{section}
\numberwithin{figure}{section}

\newcommand{\F}{\mathord{\mathbb F}}
\renewcommand{\P}{\mathord{\mathbb  P}}
\newcommand{\Q}{\mathord{\mathbb  Q}}

\newcommand{\Z}{\mathord{\mathbb Z}}

\newcommand{\CCC}{\mathord{\mathcal C}}

\newcommand{\OOO}{\mathord{\mathcal O}}

\newcommand{\RRR}{\mathord{\mathcal R}}

\newcommand{\UUU}{\mathord{\mathcal U}}
\newcommand{\VVV}{\mathord{\mathcal V}}

\font\mathgot=eufm10

\newcommand{\SSSS}{\mathord{\hbox{\mathgot S}}}

\newcommand{\dual}{\sp{\vee}}
\newcommand{\inv}{\sp{-1}}
\newcommand{\inj}{\hookrightarrow}
\newcommand{\sm}{\setminus}
\newcommand{\st}{\subset}

\newcommand{\sprime}{\sp{\prime}}
\newcommand{\sptimes}{\sp{\times}}
\newcommand{\sperp}{\sp{\perp}}
\newcommand{\ppart}{\sb{(p)}}
\newcommand{\pprimepart}{\sb{(p\sprime)}}
\newcommand{\lpart}{\sb{(l)}}

\newcommand{\threepart}{\sb{(3)}}

\newcommand{\disc}{\operatorname{\rm disc}\nolimits}
\newcommand{\Disc}[1]{\operatorname{\rm Disc}\nolimits (#1)}
\newcommand{\Discform}[1]{q\sb{#1}}
\newcommand{\Discformb}[1]{b\sb{#1}}
\renewcommand{\Im}{\operatorname{\rm Im}\nolimits}
\newcommand{\Ker}{\operatorname{\rm Ker}\nolimits}
\newcommand{\Hom}{\operatorname{\rm Hom}\nolimits}
\newcommand{\NS}{\operatorname{{\it NS}}\nolimits}
\newcommand{\PGL}{\operatorname{{\it PGL}}\nolimits}
\newcommand{\pr}{\operatorname{\rm pr}\nolimits}
\newcommand{\rank}{\operatorname{\rm rank}\nolimits}
\newcommand{\Roots}{\operatorname{{\rm Roots}}\nolimits}
\newcommand{\Sing}{\operatorname{\rm Sing}\nolimits}

\newcommand{\we}{\operatorname{\rm we}\nolimits}
\newcommand{\wt}{\operatorname{\rm wt}\nolimits}

\newcommand{\vx}{\mathord{\bf x}}
\newcommand{\vz}{\mathord{\bf 0}}

\newcommand{\set}[2]{\{\; {#1} \; \mid \; {#2} \;  \}}
\newcommand{\shortset}[2]{\{ {#1}  \mid  {#2}  \}}

\newcommand{\isom}{\smash{\mathop{\;\to\;}\limits\sp{\sim\,}}}
\newcommand{\isomto}{\smash{\mathop{\;\to\;}\limits\sp{\sim\;}}}

\newcommand{\At}{A\sb 2}
\newcommand{\Pl}{\P\sp 1}

\newcommand{\Pf}{\P\sp 4}
\newcommand{\Pll}{\P\sp 1\times \P\sp 1}
\newcommand{\tQ}{\widetilde Q}

\newcommand{\wdc}[1]{\langle #1 \rangle}

\newcommand{\nfrac}[2]{(#1)/#2}
\newcommand{\polX}{{(X, L)}}
\newcommand{\two}{\mathord{\rm II}}
\newcommand{\four}{\mathord{\rm IV}}
\newcommand{\fourstar}{\four\sp *}

\newcommand{\moduli}{\mathord{\hbox{\mathgot M}}}

\begin{document}

\title[$K3$ surfaces with ten  cusps]{$K3$ surfaces with ten cusps}

\author{Ichiro Shimada}
\address{
Department of Mathematics,
Faculty of Science,
Hokkaido University,
Sapporo 060-0810,
JAPAN
}
\email{shimada@math.sci.hokudai.ac.jp
}

\author{De-Qi Zhang}
\address{
Department of Mathematics,
National University of Singapore,
2 Science Drive 2,
Singapore 117543,
SINGAPORE
}
\email{matzdq@math.nus.edu.sg}


\begin{abstract}
We show that normal $K3$ surfaces with ten cusps exist in and only in characteristic $3$.
We determine these $K3$ surfaces according to the degrees of the polarizations.
Explicit examples are given.
\end{abstract}

\maketitle
\section{Introduction}\label{sec:introduction}
We work over an algebraically closed field $k$.

An isolated singular point of an algebraic surface
is called a \emph{cusp} if it is a rational double point of type $\At$
(Artin~\cite{Artin62, Artin66, Artin77}).

In characteristic $0$,
the  number of cusps on a normal $K3$ surface is at most nine.
Barth showed in~\cite{Barth98} that
a complex normal $K3$ surface $Y$  has nine cusps as its only singularities
if and only if 
$Y$ is  the quotient of an abelian surface by a cyclic group of order $3$.
This is a generalization of the result of~\cite{Nikulin75},
in which  Nikulin showed that 
a complex normal $K3$ surface $Y$  has sixteen nodes  as its only singularities
if and only if 
$Y$ is  the quotient of an abelian surface by the involution.
In~\cite{Barth2000},
Barth classified normal  $K3$ surfaces with nine cusps according to the degrees of the polarizations.

In positive characteristics, however,
there exist normal $K3$ surfaces $Y$ such that the singular locus $\Sing Y$ 
of $Y$ consists of \emph{ten} cusps.
The purpose of this paper is to investigate such $K3$ surfaces.

\par
\medskip 

A smooth $K3$ surface $X$
is called \emph{supersingular} (in the sense of Shioda~\cite{Shioda})  if
the N\'eron-Severi lattice $\NS (X)$ of $X$ is of rank $22$.
Supersingular $K3$ surfaces exist only in positive characteristics.
Let $X$ be a supersingular $K3$ surface in characteristic $p>0$.
Artin~\cite{Artin74} showed that
there exists a positive integer $\sigma (X)\le 10$
such that
$\disc \NS (X)=-p\sp{2\sigma (X)}$ holds.
This integer $\sigma (X)$ is  called the \emph{Artin invariant of $X$}.

We denote by
 $U(m)$ 
the lattice of rank $2$
whose intersection matrix is equal to 
$$
\begin{pmatrix}
0 & m \\
m & 0
\end{pmatrix}.
$$

Our main results are Theorems 1.1 and  1.4 - 1.6.

\begin{theorem}\label{thm:A}
Let $Y$ be a normal $K3$ surface such that $\Sing Y$ consists of ten cusps,
and  $\rho: X\to Y$ the minimal resolution of $Y$.
Let  $R\sb{\rho}$ be the sublattice of $\NS (X)$  generated by the classes of
the $(-2)$-curves that are contracted by $\rho$.
Then the following hold:

{\rm (1)}
The characteristic of the ground field $k$ is $3$.

{\rm (2)}
The orthogonal complement $R\sb\rho\sperp$ of  $R\sb{\rho}$ in $\NS (X)$
is isomorphic to  either $U(1)$ or $U(3)$.

{\rm (3)}
If $R\sb\rho\sperp\cong U(1)$, then $\sigma (X)\le 5$,
while 
if $R\sb\rho\sperp\cong U(3)$, then $\sigma (X)\le 6$.
\end{theorem}

Before we state the other main results, we fix the terminology below.

\begin{definition}
Let $L$ be a line bundle on a smooth $K3$ surface $X$.
We say that $L$ is \emph{very ample modulo $(-2)$-curves}
if the following hold:
\begin{itemize}
\item[{\rm (i)}] The complete linear system $|L|$ has no fixed components,
and hence has no base points by~\cite[Corollary 3.2]{SD}.
In particular, $|L|$ defines a morphism $\Phi\sb{|L|} : X\to \P\sp{N}$, where $N=L^2/2+1$.
\item[{\rm (ii)}] The morphism $\Phi\sb{|L|}$ is birational onto the image
$Y\sb{(X, L)}:=\Phi\sb{|L|} (X)$.
\end{itemize}
 A \emph{polarized  $K3$ surface}  is a pair $(X, L)$ of a $K3$ surface $X$ and a line bundle $L$
on $X$ that  is very ample modulo $(-2)$-curves.
The \emph{degree} of a polarized  $K3$ surface $(X, L)$  is defined to be $L^2$.
\end{definition}
\begin{definition}
Let $(X, L)$ be a polarized  $K3$ surface.
We denote by 
$$
\rho\sb L: X\to Y\sb{\polX}
$$
the birational morphism induced by $|L|$.
By~\cite[Theorem 6.1]{SD},
$\rho\sb L$ is a contraction of an $ADE$-configuration of $(-2)$-curves on $X$.
We denote by  $R\sb{(X, L)}$ the sublattice of $\NS (X)$ generated by the classes of
the $(-2)$-curves that are contracted by $\rho\sb L$.
We also denote by $\RRR\sb{(X, L)}$ the $ADE$-type of the configuration of 
these $(-2)$-curves.
\end{definition}
Note that $\RRR\sb{(X, L)}=10\At$ is equivalent to saying that $\Sing Y\sb{\polX}$ consists of ten cusps. The degree of $(X, L)$ can be completely determined:
\begin{theorem}\label{thm:U1}
The following conditions on a positive integer  $d$ are equivalent:
\begin{itemize}
\item[{\rm (i)}]
$d=2ab$, where $a$ and $b$ are integers $\ge 3$ such that $a\ne b$.
 \item[{\rm (ii)}]
There exists a polarized supersingular $K3$ surface $(X, L)$
of degree $d$
such that  $\RRR\sb{(X, L)}=10\At$ and 
$R\sb{\polX}\sperp\cong U(1)$.
\item[{\rm (iii)}]
Every supersingular $K3$ surface $X$ in characteristic $3$ with $\sigma (X)\le 5$
admits a line bundle  $L$  such that 
$(X, L)$ is a polarized $K3$ surface of degree $d$ satisfying 
$\RRR\sb{(X, L)}=10\At$ and 
$R\sb{\polX}\sperp\cong U(1)$.
\end{itemize}
\end{theorem}
\begin{theorem}\label{thm:U3}
The following conditions on a positive integer  $d$ are equivalent:
\begin{itemize}
\item[{\rm (i)}]
$d\equiv 0 \bmod 6$.
 \item[{\rm (ii)}]
There exists a polarized supersingular $K3$ surface $(X, L)$
of degree $d$
such that  $\RRR\sb{(X, L)}=10\At$ and 
$R\sb{\polX}\sperp\cong U(3)$.
\item[{\rm (iii)}]
Every supersingular $K3$ surface $X$ in characteristic $3$ with $\sigma (X)\le 6$
admits a line bundle  $L$  such that 
$(X, L)$ is a polarized $K3$ surface of degree $d$ satisfying 
$\RRR\sb{(X, L)}=10\At$ and 
$R\sb{\polX}\sperp\cong U(3)$.
\end{itemize}
\end{theorem}
Supersingular $K3$ surfaces with ten cusps can be obtained as purely inseparable triple covers of 
the smooth quadric surface $Q=\Pll$.
From now on to the end of this paragraph,
we assume that $k$ is of characteristic $3$.
For integers $a$ and $b$,
we denote by $\OOO\sb{Q} (a, b)$ the invertible sheaf 
$\pr\sb 1 \sp * \OOO\sb{\Pl} (a) \otimes \pr\sb 2 \sp * \OOO\sb{\Pl} (b)$ of $Q=\Pll$,
and by $L\sb{Q} (a, b)\to Q=\Pll$ the corresponding line bundle.
Because we are in characteristic $3$,
the differential map
$$
d : H\sp 0 (Q, \OOO\sb{Q} (3, 3)) \to H\sp 0 ( Q, \Omega\sb{Q} \sp 1 (3, 3))
$$
is well-defined by the isomorphism $L\sb{Q} (3, 3)\cong L\sb{Q}(1, 1)\sp{\otimes 3}$.
For $G\in  H\sp 0 (Q, \OOO\sb{Q} (3, 3))$, we denote by $Z(dG)$ the subscheme of $Q$ defined by $dG=0$.
If  $\dim Z(dG)=0$, then
$$
\textrm{length}\; \OOO\sb{Z(dG)}=c\sb 2 ( \Omega\sb{Q} \sp 1 (3, 3)) =10
$$
holds, where $c\sb 2$ is the second Chern class.
We put
$$
\UUU\sb{3,3}:=\set{G\in  H\sp 0 (Q, \OOO\sb{Q} (3, 3)) }{\textrm{$Z(dG)$ is reduced and of dimension $0$}},
$$
which is a Zariski open dense subset of $H\sp 0 (Q, \OOO\sb{Q} (3, 3))$.
For a non-zero   $G\in H\sp 0 (Q, \OOO\sb{Q} (3, 3))$, we denote by
$$
\pi\sb G : Y\sb G \to Q=\Pll
$$
the purely inseparable triple cover of $Q$ defined by 
$$
W^3=G,
$$
where $W$ is a fiber coordinate of the line bundle $L\sb{Q}(1, 1)$.
It is easy to see that $G$ is contained in $\UUU\sb{3,3}$ if and only if 
$Y\sb G$ is a normal $K3$ surface
such that $\Sing  Y\sb G =\pi\sb G\inv (Z(dG))$
consists of ten cusps.
In particular, if $G\in \UUU\sb{3,3}$, then  
the minimal resolution $X\sb G$ of $Y\sb G$ is a supersingular $K3$ surface
with Artin invariant $\le 6$ by Theorem~\ref{thm:A}.
Conversely, we have the following:
\begin{theorem}\label{thm:insep}
Let $X$ be a supersingular $K3$ surface in characteristic $3$ with Artin invariant $\le 6$.
Then there exists $G\in \UUU\sb{3,3}$ such that $X$ is isomorphic to $X\sb G$.
\end{theorem}
\noindent
We put
$$
 \VVV\sb{1,1}:=\set{H^3\in H^0 (Q, \OOO\sb{Q}(3,3))}{H\in H^0(Q, \OOO\sb{Q}(1,1))},
$$
which is an additive group 
acting  on $\UUU\sb{3,3}$ by $G\mapsto G+H^3$ $(G\in \UUU\sb{3,3}, H^3\in \VVV\sb{1,1})$.
For $G, G\sprime\in \UUU\sb{3, 3}$, the triple covers
$Y\sb G$ and $Y\sb{G\sprime}$ are isomorphic over $Q$ if and only if
$G=cG\sprime +H^3$ holds for some $c\in k\sptimes$ and $H^3\in \VVV\sb{1,1}$.
Hence the space 
$$
\moduli:=(\PGL (2, k)\times \PGL(2, k))\backslash \P\sb * (\UUU\sb{3,3}/ \VVV\sb{1,1})
$$
is a moduli space of supersingular $K3$ surfaces in characteristic $3$
with Artin invariant $\le 6$.
We remark that, 
since $\dim \UUU\sb{3,3}=16$ and  $\dim \VVV\sb{1,1}=4$,
we have 
$$
\dim \moduli = 16-4-1-(3+3)=5,
$$
as is predicted from the result of Artin~\cite{Artin74}.
In particular, the {\it unique} supersingular $K3$ surface of Artin
invariant $1$ has the following precise model:

\begin{example}\label{example:G0}
We put
$$
G\sb 0:=(x^3-x)(y^3-y),
$$
where $x$ and $y$ are affine coordinates of the two factors of $Q=\Pll$.
Then $Z(dG\sb 0)$ is equal to 
$$ 
\set{(\alpha, \beta)}{\alpha, \beta\in \F\sb 3 }\cup \{(\infty, \infty)\}.
$$
Therefore $G\sb 0\in \UUU\sb{3,3}$.
It can be shown that
the Artin invariant  of the supersingular $K3$ surface $X\sb{G\sb 0}$
is  $1$.
See  Example~\ref{example:sigmaG0}.
\end{example}
\par
\medskip
Supersingular $K3$ surfaces in characteristic $2$ with $21$ nodes are investigated in
\cite{Shimada2003, Shimada2004, Shimada2004moduli}.
In particular,
it was shown there that every supersingular $K3$ surface in characteristic $2$
is birational to a purely inseparable double cover of the projective plane
with $21$ nodes;
that is, 
every supersingular $K3$ surface in characteristic $2$ is obtained as a generic Zariski surface~\cite{BL}.

Quasi-elliptic $K3$ surfaces in characteristic $3$  with a section and ten singular fibers of type 
$\At^*$
are constructed explicitly in~\cite{Ito92}.
The Artin invariants of these supersingular $K3$ surfaces are $\le 5$.

A family of smooth quartic surfaces in characteristic $3$  containing an $ADE$-configuration
of lines of type $10\At$ is constructed in~\cite{Shimada92}.
A general member of the family is of Artin invariant $6$.
See Example~\ref{example:quartic} for details.
\par
\medskip
This paper is organized as follows.
In \S\ref{sec:disc},
we review the theory of discriminant forms of lattices due to Nikulin~\cite{Nikulin79}.
In \S\ref{sec:NS},
we quote from Artin~\cite{Artin74}, Rudakov-Shafarevich~\cite{RS},
 Saint-Donat~\cite{SD} and Nikulin~\cite{Nikulin91}
some known facts about N\'eron-Severi lattices and 
polarizations of supersingular $K3$ surfaces.
In \S\ref{sec:thmA},
we prove Theorem~\ref{thm:A} using the theory of discriminant forms.
In \S\ref{sec:U3} and \S\ref{sec:U1},
we prove Theorems~\ref{thm:U3}~and~\ref{thm:U1}.
We reduce the problem of  existence of the polarizations on a supersingular $K3$ surface
to a problem of  existence of ternary codes with  certain properties,
and solve the latter  by computer.
In \S\ref{sec:insep},
we prove Theorem~\ref{thm:insep}.
The proof presented here seems to be quite lattice-intensive.
We think there should be a more elementary proof. 
See Question~\ref{question1}.
\par
\medskip
{\bf Acknowledgment.} This work was done during the second author's visit
to Hokkaido University who likes to express his thanks for the very warm hospitality.
\section{Discriminant forms of lattices}\label{sec:disc}
For a finite abelian group $A$ and a prime integer $p$,
we denote by
$$
A=A\ppart\times A\pprimepart
$$
the decomposition of $A$ into the $p$-part $A\ppart$ and the  $p$-prime-part $A\pprimepart$ of $A$.

A lattice is, by definition,
a free $\Z$-module of finite rank with a non-degenerate symmetric 
$\Z$-valued bilinear form.
A lattice $\Lambda$ is said to be \emph{even} if $v^2\in 2\Z$ holds for every $v\in \Lambda$.
Let $\Lambda$ be an even lattice.
We denote by $\Lambda\dual$ the \emph{dual lattice} $\Hom (\Lambda, \Z)$.
We have a natural embedding $\Lambda\inj \Lambda\dual$
of finite cokernel, and a symmetric bilinear form
$\Lambda\dual\times\Lambda\dual\to\Q$ that extends the $\Z$-valued symmetric bilinear form on $\Lambda$.
We put
$$
\Disc {\Lambda} :=\Lambda\dual/\Lambda,
$$
and call it the \emph{discriminant group of $\Lambda$}.
We then define the \emph{discriminant form}
\begin{eqnarray*}
\Discform{\Lambda} &:& \Disc{\Lambda} \to \Q/2\Z \quad\textrm{and}\\
\Discformb{\Lambda} &:& \Disc{\Lambda} \times \Disc{\Lambda} \to \Q/\Z
\end{eqnarray*}
by
\begin{eqnarray*}
\Discform{\Lambda} (\bar v)&:=&v^2\bmod 2\Z\quad\textrm{and}\\
\Discformb{\Lambda} (\bar v, \bar w)&:=&vw\bmod \Z=
(\Discform{\Lambda}(\bar v +\bar w) -\Discform{\Lambda}(\bar v)-\Discform{\Lambda}(\bar w))/2,
\end{eqnarray*}
where $v, w\in\Lambda\dual$, and $\bar v:=v \bmod \Lambda$, $\bar w :=w \bmod \Lambda$.
Let $p$ be a prime integer dividing $| \Disc{\Lambda}| =|\disc\Lambda|$.
Then $\Disc{\Lambda}\ppart$ and $\Disc{\Lambda}\pprimepart$ are
orthogonal with respect to $\Discformb{\Lambda}$.
We put
\begin{eqnarray*}
{\Discform{\Lambda}}\ppart:= \Discform{\Lambda}| \Disc{\Lambda}\ppart, &&
{\Discform{\Lambda}}\pprimepart:= \Discform{\Lambda}| \Disc{\Lambda}\pprimepart, \\
{\Discformb{\Lambda}}\ppart:= \Discformb{\Lambda}| \Disc{\Lambda}\ppart\times \Disc{\Lambda}\ppart, &&
{\Discformb{\Lambda}}\pprimepart:= \Discformb{\Lambda}| \Disc{\Lambda}\pprimepart\times \Disc{\Lambda}\pprimepart.
\end{eqnarray*}
For a subgroup $H$ of $\Disc{\Lambda}$,
we denote by $H\sperp$ the orthogonal complement of $H$ with respect to $\Discformb{\Lambda}$.
Note that $(H\sperp)\ppart$ is canonically isomorphic to
$$
(H\ppart)\sperp:=\set{x \in \Disc{\Lambda}\ppart}{\hbox{${\Discformb{\Lambda}}\ppart (x, y)=0$ for any $y\in H\ppart$ }}.
$$
We  will use the notation $H\sperp\ppart$
to denote $(H\sperp)\ppart=(H\ppart)\sperp$.
A subgroup $H\subset \Disc{\Lambda}$ is called \emph{isotropic} if $\Discform{\Lambda}|H$ is constantly equal to $0$.
If $H$ is isotropic, then $H$ is contained in $H\sperp$.
Note that we have
$$
(H\sperp/H)\ppart=H\sperp\ppart/H\sp{\phantom{\perp}}\ppart.
$$
An \emph{overlattice} of $\Lambda$ is, by definition,
a submodule $\Lambda\sprime$ of $\Lambda\dual$ containing $\Lambda$ such that the $\Q$-valued symmetric bilinear form on
$\Lambda\dual$ takes values in $\Z$ on $\Lambda\sprime$.
\begin{proposition}[Nikulin \cite{Nikulin79}]\label{prop:nikulin}
Let
$\pr\sb{\Lambda} : \Lambda\dual \to \Disc{\Lambda}$
be the natural projection.
The correspondence
$$
H\mapsto \Lambda\sb H :=\pr\sb{\Lambda}\inv (H)
$$
gives a bijection from 
the set of isotropic subgroups of $\Disc{\Lambda}$
to the set of even overlattices of $\Lambda$.
For an isotropic subgroup $H$, 
the discriminant group of $\Lambda\sb{H}$ is isomorphic to $H\sperp/H$.
\end{proposition}
\begin{remark}
If $\Lambda$ is of rank $r$, then $\Disc{\Lambda}$ is generated by $r$ elements.
\end{remark}
%
%
%
%
\par
\medskip
A vector $v$ in an even  \emph{negative-definite} lattice $\Lambda$ is called a \emph{root}
if $v^2=-2$.
We denote by
$\Roots (\Lambda)$ the set of roots in $\Lambda$.
It is known that $\Roots (\Lambda)$
forms a root system of type $ADE$~(\cite{B, E}).
An even  negative-definite lattice $\Lambda$ is called a \emph{root lattice} if it is generated by $\Roots (\Lambda)$.
\par
\medskip
Let $\Z [10 \At ]$ denote the  root lattice of type $10\At$.
Then $\Z [10 \At ]$ is generated by roots
$c\sb i, d\sb i$ $(i=1, \dots, 10)$
satisfying
$$
c\sb i^2=d\sb i^2=-2, \quad c\sb i d\sb i=1, \quad\textrm{and}\quad
 \langle c\sb i, d\sb i\rangle \perp \langle c\sb j, d\sb j\rangle \;\;\textrm{if $i\ne j$}.
$$
We have 
$$
\Roots(\Z[10\At])=\{ \pm c\sb i, \pm d\sb i, \pm (c\sb i + d\sb i)\;\;\; (i=1, \dots, 10)\},
$$
and
$$
\Z[10\At]\dual=\set{\sum\sb{i=1}\sp{10} (s\sb i c\sb i+t\sb i d\sb i)/3}{s\sb i, t\sb i\in \Z, \; 
s\sb i+t\sb i\equiv 0 \bmod 3\;\; (i=1, \dots, 10)}.
$$
We put
$$
\gamma\sb i := (c\sb i + 2 d\sb i)/3\;\bmod\; \Z [10 \At ]\;\;\in\;\; \Disc{\Z[10\At]}.
$$
Then we have 
$$
\Disc{\Z[10\At]}=\F\sb 3 \gamma\sb 1 \oplus \cdots \oplus \F\sb 3 \gamma\sb{10},
$$
and
\begin{equation}\label{eq:discform10At}
\Discform{\Z[10\At]}(x\sb 1 \gamma\sb 1+\cdots+x\sb{10}\gamma\sb{10}) =
-2 (x\sb 1^2 + \cdots + x\sb{10}^2)/3 \;\;\in\;\;\Q/2\Z.
\end{equation}
For a vector 
$$
\vx=(x\sb 1, \dots, x\sb{10})=x\sb 1 \gamma\sb 1+\cdots+x\sb{10}\gamma\sb{10}\;\in\;  \Disc{\Z[10\At]}\cong \F\sb 3\sp{10},
$$ 
we define the \emph{Hamming weight} $\wt (\vx)$  of $\vx$ by
$$
\wt (\vx ):=|\set{i}{x\sb i\ne 0}|\;\;\in\;\;\Z\sb{\ge 0}.
$$
Then, for a vector $r\in \Z[10\At]\dual$, we have
\begin{equation}\label{eq:wtr}
r^2\le -(2/3) \wt(\bar r), \quad\textrm{where $\bar r :=r \bmod \Z[10\At] \;\in\; \Disc{\Z[10\At]}$}.
\end{equation}
Moreover, 
\begin{equation}\label{eq:wtrinv}
\parbox{10cm}{
for a vector $\vx\in  \Disc{\Z[10\At]}$,
there exists a vector $r\in \Z[10\At]\dual$ such that
$\bar r =\vx$ and 
$r^2= (-2/3) \wt(\vx)$ hold.
}
\end{equation}
\par
\bigskip
Let $e$ and $f$ be basis of the lattice $U(m)$ satisfying
$$
e^2=f^2=0, \qquad ef=m.
$$
We put $e\dual:= f/m$ and $f\dual :=e/m$.
Then $\Disc{U(m)}\cong (\Z/ m\Z)\sp 2 $ is generated by 
$$
\bar e\dual:=e\dual \bmod U(m) \quad\textrm{and}\quad 
\bar f\dual:=f\dual \bmod U(m),
$$
and the discriminant form  is given by
\begin{equation}\label{eq:discformUm}
\Discform{U(m)}(y\sb 1 \bar e\dual + y\sb 2  \bar f\dual)= 2 y\sb 1 y\sb 2 /m \;\;\in\;\;\Q/2\Z.
\end{equation}
\section{N\'eron-Severi lattices of supersingular $K3$ surfaces}\label{sec:NS}
A lattice $\Lambda$ is called \emph{hyperbolic}
if the signature of $\Lambda$ is $(1, \rank \Lambda -1)$.
Let $p$ be a prime integer.
A lattice $\Lambda$ is called \emph{$p$-elementary}
if $\Disc{\Lambda}$ is a $p$-elementary abelian group;
that is, $p \Lambda\dual \subseteq \Lambda$ holds.
An overlattice of a hyperbolic $p$-elementary lattice is 
again  hyperbolic and $p$-elementary.
\par
\medskip
The following is due to Artin~\cite{Artin74} and Rudakov-Shafarevich~\cite{RS}.
\begin{theorem}\label{thm:ARS}
Let $X$ be a supersingular $K3$ surface in characteristic $p>0$.
Then $\NS (X)$ is an even  hyperbolic $p$-elementary lattice.
\end{theorem}
The following is due to Rudakov-Shafarevich~\cite[Section 1]{RS}.
\begin{theorem}\label{thm:NRS}
Suppose that $p$ is odd.
Let $\sigma$ be a positive integer $\le 10$.
Then the lattice $N$ with the following properties is unique up to isomorphisms:

{\rm (i)} $N$ is even, hyperbolic of rank $22$, and

{\rm (ii)} $\Disc{N} \cong \F\sb p \sp{2\sigma}$.
\end{theorem}
From now on to the end of this section,
we assume that $p$ is \emph{odd}.
We denote the lattice $N$ in Theorem~\ref{thm:NRS} by $N\sb{p, \sigma}$.
Let $X$ be a supersingular $K3$ surface in  characteristic $p$
with $\sigma (X)=\sigma$.
By  Theorems~\ref{thm:ARS} and~\ref{thm:NRS},
there exists an isometry
$$
\phi : N\sb{p, \sigma} \isomto \NS(X).
$$
More precisely, we have the following:
\begin{proposition}\label{prop:h}
Let $h$ be a vector of $N\sb{p, \sigma}$ such that $h^2\ge 4$,
and let $X$ be a supersingular $K3$ surface in  characteristic $p$
with $\sigma (X)=\sigma$.

{\rm (1)}
The following conditions are equivalent:
\begin{itemize}
\item[{\rm (i)}]
There exist no vectors $u\in N\sb{p, \sigma}$ satisfying  $hu=1$ or $2$ and $u^2=0$, 
and there exist no vectors $b\in N\sb{p, \sigma}$ satisfying  $h=2b$ and $b^2=2$.
\item[{\rm (ii)}]
There exists an isometry $\phi : N\sb{p, \sigma} \isomto  \NS(X)$
such that $\phi (h)$ is the class $[L]$ of a line bundle $L$ that is very ample
modulo $(-2)$-curves.
\end{itemize}

{\rm (2)}
Suppose that the conditions in {\rm (1)} are fulfilled,
and let $L$ be a line bundle very ample
modulo $(-2)$-curves such that $\phi (h)=[L]$ by some isometry $\phi$.
Then $Y\sb{(X, L)}$ has only rational double points as its singularities,
and the $ADE$-type $\RRR\sb{(X, L)}$ of 
$\Sing Y\sb{\polX}$ is equal to that of the root system
$$
\Roots (h\sperp):=\set{r\in N\sb{p, \sigma}}{rh=0,\; r^2=-2}.
$$
\end{proposition}
For the proof, we use the following results due to 
Nikulin~\cite[Proposition 0.1]{Nikulin91} and  
Saint-Donat~\cite[Section 5]{SD}.
\begin{proposition}[Nikulin \cite{Nikulin91}]\label{prop:nikulin1}
Let $L$ be a nef line bundle on a $K3$ surface $X$
with $L^2>0$.
If $|L|$ has a  fixed component,
then 
$|L|$ is equal to  $m|U|+\Gamma$,
where $\Gamma$ is the fixed part of $|L|$,
$|U|$ is an elliptic pencil, 
and 
$U^2=0$, $U\Gamma=1$, $\Gamma^2=-2$, $m=\dim |L|=L^2/2+1$ hold.
If $|L|$ has no fixed components,
then a general member of $|L|$ is irreducible and $\dim |L|=L^2/2+1$.
\end{proposition}
\begin{proposition}[Saint-Donat \cite{SD}]\label{prop:SD1}
Let $|L|$ be a complete linear system without fixed components  on a $K3$ surface $X$
such that  $L^2\ge 4$.
Then the morphism $\Phi\sb{|L|}$ fails to be birational onto its image 
if and only if one of the following holds:
\begin{itemize}
\item[{\rm (i)}]
There exists an irreducible curve $U$ such that $U^2=0$ and $UL=2$.
\item[{\rm (ii)}]
There exists an irreducible curve $B$ such that $B^2=2$ and $L=\OOO\sb X (2B)$.
\end{itemize} 
\end{proposition}
\begin{proof}[Proof of Proposition~\ref{prop:h}]
The assertion  (2) follows from~\cite[Theorem 6.1]{SD}
and~\cite[Lemma 2.4]{Shimada2003}. We now prove (1).

\par
Suppose that the condition (i) in (1) holds.
By~\cite[Section 3, Proposition 3]{RS},
there exists an isometry $\phi : N\sb{p, \sigma} \isom \NS (X)$ such that $\phi (h)$ is 
the class of a \emph{nef} line bundle $L$.
By Proposition~\ref{prop:nikulin1},
$|L|$ is fixed component free.
By Proposition~\ref{prop:SD1},
$\Phi \sb{|L|}$ is birational onto its image. So (ii) is true.

\par
Conversely, suppose that (ii) holds.
We assume that there exists a vector $u\in N\sb{p, \sigma}$ satisfying $hu=1$ or $2$ and $u^2=0$,
and derive a contradiction by the argument in~\cite[Proof of Proposition 1.7]{U}.
By the Riemann-Roch theorem,
$\phi (u)$ is the class $[U]$ of an effective divisor $U$
such that $\dim |U|\ge 1$.
Let $D+\Delta$ be a general member of $|U|$,
where $\Delta$ is the fixed part of $|U|$.
We have $D\ne 0$ and $D^2\ge 0$.
If $DL=0$, then $D^2<0$ would follow by Hodge index theorem,
a contradiction.
Since $L$ is nef, $\Delta L\ge 0$.
Therefore, we have $DL=1$ or $2$.
Then the image of $D$ by $\Phi\sb{|L|}$ is either a line or a plane conic.
In any case, we have $\dim |D|=0$,
which is a contradiction.

\par
Next we assume that there exists a vector $b\in N\sb{p, \sigma}$
such that $h=2b$ and $b^2=2$.
Let $B$ be an effective divisor such that $\phi (b)=[B]$.
Since $[B]=[L]/2$,
$B$ is nef.
If there exists an irreducible member in $|B|$,
then Proposition~\ref{prop:SD1} implies that $\Phi\sb{|L|}$ is not birational 
onto its image.
If there exist no irreducible members in $|B|$,
then  Proposition~\ref{prop:nikulin1} implies that 
$|B|$ has a fixed component, and $|B|$ is written as $2 |U|+\Gamma$,
where $UB=1$ and $U^2=0$.
Then $UL=2$ follows. Hence $\Phi\sb{|L|}$ is not birational onto its image,
and we get a contradiction. So (i) is true.
Thus the assertion (1) is proved.
\end{proof}
\begin{remark}\label{rem:degree8}
If there exists a vector $b$ such that $h=2b$ and $b^2=2$,
then $h$ is of degree $8$.
\end{remark}
\section{Proof of Theorem~\ref{thm:A}}\label{sec:thmA}
Theorem~\ref{thm:A} follows from the structure theorem  of N\'eron-Severi lattices of supersingular $K3$ surfaces
(Theorems~\ref{thm:ARS} and~\ref{thm:NRS}), and a purely lattice-theoretic Lemma~\ref{lemma:lattice} below.
A sublattice $\Lambda\sprime\st \Lambda$ is called \emph{primitive in $\Lambda$} if
$(\Lambda\sprime\otimes\Q) \cap \Lambda=\Lambda\sprime$ holds.
\begin{lemma}\label{lemma:lattice}
Let $N$ be an even hyperbolic $p$-elementary lattice of rank $22$
such that $\Disc{N}$ is isomorphic to $\F\sb p \sp{ 2 \sigma}$,
where $\sigma$ is a positive integer.
Suppose that $N$ contains a sublattice $R$ isomorphic to $\Z [10\At]$.
Then $p=3$,  and the orthogonal complement $R\sperp$ of $R$ in $N$ is isomorphic to $U(1)$ or $U(3)$.
If $S\cong U(1)$, then $\sigma\le 5$,
while if $S\cong U(3)$, then $\sigma\le 6$.
\end{lemma}
\begin{proof}
We put $S:=R\sperp$,
which is an even hyperbolic lattice of rank $2$  primitive in $N$.
Then $N$ is an overlattice of the orthogonal direct sum $R\oplus S$.
We put
$$
H:= N/(R\oplus S).
$$
Clearly, we may assume that $H \ne (0)$.

\par
Note that $H$ is an isotropic subgroup of $ \Disc{R\oplus S}=\Disc{R} \oplus \Disc{S}$ with respect to
$\Discform {R\oplus S}=\Discform{R}\oplus\Discform{S}$,
and 
  $\Disc{N}\cong H\sperp/H$ is  a $p$-elementary abelian group.
Since $S$ is primitive in $N$,
we have
\begin{equation}\label{eq:HcapS}
H\cap (0 \oplus \Disc{S}) =0.
\end{equation}
Let $l$ be a prime integer different from $3$ and $p$.
Assume that $\Disc{S}\lpart$ is not $0$.
Since $\Disc{N}\lpart=0$,
we see that  $H\lpart$
is not $0$.
Since $\Disc {\Z[10\At]}\lpart=0$,
we have $H\lpart\subset (0 \oplus \Disc{S}\lpart)$,
which contradicts~\eqref{eq:HcapS}.
Hence we obtain
\begin{equation}\label{eq:DiscSlpart}
\Disc{S}\lpart=0 \quad\textrm{for any prime $l$ distinct from $ 3$ and $ p$}.
\end{equation}
Let $m\sb 3 : \Disc{S}\threepart \to \Disc{S}\threepart$
be the homomorphism given by $m\sb 3 (x):=3x$.
Since every element of $\Disc{R}$ is annihilated by multiplication by $3$,
the image  $H\threepart\sp S\subset  \Disc{S}\threepart$ of 
 $H\threepart\subset \Disc{R}\threepart\oplus \Disc{S}\threepart$ 
by the projection to the factor $\Disc {S}\threepart$
is contained in $\Ker m\sb 3$ by~\eqref{eq:HcapS}:
\begin{equation}\label{eq:inKer}
H\threepart\sp S \;\;\subseteq\;\; \Ker m\sb 3.
\end{equation}
Therefore,  $\Im m\sb 3$ is contained in the orthogonal complement
of  $H\threepart\sp S$
with respect to ${\Discform{S}}\threepart$.
Hence we obtain
\begin{equation}\label{eq:Imm}
0\oplus \Im m\sb 3 \subset H\sperp\threepart.
\end{equation}
\par
\medskip

We assume $p\ne 3$, and derive a contradiction.
By~\eqref{eq:DiscSlpart},
we have
\begin{equation}\label{eq:prod}
\Disc{S} = \Disc {S}\threepart \times \Disc{S}\ppart.
\end{equation}
Since $\Disc{R}\ppart=0$,
the property~\eqref{eq:HcapS} implies $H\ppart=0$.
Therefore $\Disc{N}=\Disc{N}\ppart$ is isomorphic to $\Disc{S}\ppart$.
Since $\dim \sb{\F\sb p} \Disc{N}=2\sigma$ is positive and even, 
and $S$ is of rank $2$,
we obtain
\begin{equation}\label{eq:ASp}
\Disc {S}\ppart \cong \F\sb p\sp{ 2}.
\end{equation}
On the other hand,
from $\Disc{N}\threepart=0$, we obtain
\begin{equation}\label{eq:H3}
H\threepart=H\sperp\threepart.
\end{equation}
By~\eqref{eq:HcapS},~\eqref{eq:Imm} and~\eqref{eq:H3},
we obtain $\Im m\sb 3 =0$; that is, $\Disc{S}\threepart$ is $3$-elementary.
From~\eqref{eq:H3},
we have 
$10 + \dim \sb{\F\sb 3} \Disc{S}\threepart=2 \dim \sb{\F\sb 3} H\threepart$,
and hence
 $\dim\sb{\F\sb 3 } \Disc{S}\threepart$ is even.
Since $S$ is of rank $2$,
we obtain  
\begin{equation}\label{eq:AS3}
\Disc {S}\threepart \cong 0 \;\;\textrm{or}\;\; \F\sb 3\sp{ 2}.
\end{equation}

Suppose that $\Disc {S}\threepart \cong 0$.
Then $H\threepart$ can be regarded as an isotropic subgroup
of $\Disc{R}$ with respect to $\Discform{R}$.
Because $H\threepart=H\sperp\threepart$, 
the corresponding overlattice of $R$ would be an even unimodular negative-definite lattice
of rank $20$.
This contradicts the classification of unimodular lattices~(\cite[Chapter V]{Serre}).

Suppose that $\Disc {S}\threepart \cong \F\sb 3\sp{ 2}$.
By~\eqref{eq:prod} and~\eqref{eq:ASp},
$S$ is an even indefinite lattice of rank $2$
such that $\Disc{S}\cong (\Z/3p\Z)\sp{ 2}$.
By the classification of indefinite lattices of rank $2$~(\cite[Chapter 15, Section 3]{CS}),
we see that the intersection matrix of $S$ with respect to an appropriate basis is 
$$
\begin{pmatrix}
0 & 3p \\
3p & 0
\end{pmatrix},
\qquad
\textrm{or}
\qquad
p=2\;\;\textrm{and}\;\;
\begin{pmatrix}
6 & 6 \\
6 & 0
\end{pmatrix}.
$$
In any case, the quadratic  form $(\Disc{S}\threepart, {\Discform{S}}\threepart)$ is isomorphic  to 
$$
\left(\:\F\sb 3\sp{ 2}, \: \left[ \begin{matrix} 0 & 1/3 \\ 1/3 & 0 \end{matrix} \right] \:\right)
\;\cong\;
(\Disc{U(3)}, \Discform{U(3)}).
$$
Therefore the isotropic subgroup $H\threepart$ of $\Disc {R}\oplus \Disc{S}\threepart$
satisfying  $H\threepart=H\sperp\threepart$ would yield an even hyperbolic unimodular lattice of rank $22$
as an overlattice of $R\oplus U(3)$,
which again contradicts the classification of unimodular lattices.

Therefore $p=3$ is proved.
\par
\medskip
By~\eqref{eq:DiscSlpart}, we have 
$\Disc {S}=\Disc {S}\threepart$,
and hence $H=H\threepart$ holds.
Suppose that $(\xi, \eta)\in H\sperp$,
where $\xi\in \Disc{R}$ and $\eta\in \Disc{S}$.
Since  $H\sperp/H$ is $3$-elementary, 
we have $(3\xi, 3\eta)=(0, 3\eta)\in H$.
By~\eqref{eq:HcapS},
we have $3\eta=0$.
Therefore
the image  $(H\sperp )\sp S\subset  \Disc{S}$ of 
 $H\sperp\subset \Disc{R}\oplus \Disc{S}$ 
by the projection to the factor $\Disc {S}$
is contained in $\Ker m\sb 3$:
\begin{equation}\label{eq:inKer2}
(H\sperp)\sp S \;\;\subset\;\; \Ker m\sb 3.
\end{equation}
Next we will show that $S$ is isomorphic to $U(1)$ or  $U(3)$.
Since $H\sperp/H$ is $3$-elementary,
\eqref{eq:HcapS} and~\eqref{eq:Imm} implies that $m\sb 3 (\Im m\sb 3)=0$;
that is, $9x=0$ for any $x\in \Disc {S}$.
Since 
$$
2\sigma =\dim\sb{\F\sb 3} (H\sperp/H) =10 + \log\sb 3 | \Disc{S}| - 2 \log \sb 3 |H|
$$
is even and $S$ is of rank $2$,
$\Disc {S}$ is isomorphic to $0$, $\F\sb 3\sp 2$, $\Z/9\Z$ or $(\Z/9\Z)\sp 2$.

We first assume  that $\Disc{S}$ is a cyclic group of order $9$, and derive a contradiction.
Let $\gamma$ be a generator of $\Disc{S}$.
We have $\Im m\sb 3 =\Ker m\sb 3 =\langle 3\gamma\rangle$.
Let $H\sp R\subset \Disc{R}$ and $H\sp S\subset \Disc{S}$ be the images of  
$H\subset \Disc{R}\oplus \Disc{S}$ by the projections to 
the factors $\Disc{R}$ and $\Disc{S}$, respectively.
\begin{claim}\label{claim:1}
We have
$$
H\sperp=(H\sp R)\sperp \oplus (H\sp S)\sperp,
$$
where 
 $(H\sp R)\sperp\subset \Disc{R}$ and $(H\sp S)\sperp\subset \Disc{S}$ are 
the orthogonal complements of $H\sp R$ and $H\sp S$
with respect to $\Discform{R}$ and $\Discform{S}$, respectively.
In particular, we have $ (H\sp S)\sperp=(H\sperp)\sp S$.
\end{claim}
\begin{proof}
It is obvious that 
$H\sperp$ contains $(H\sp R)\sperp \oplus (H\sp S)\sperp$.
Suppose that $(\xi, \eta)\in H\sperp$, where $\xi\in \Disc{R}$ and $\eta\in \Disc{S}$.
By~\eqref{eq:inKer2}, we have $\eta\in\Ker m\sb 3=\Im m\sb 3$.
By~\eqref{eq:Imm}, we have  $(0, \eta)\in H\sperp$ and hence $(\xi, 0)\in H\sperp$
hold.
Because $(\xi, 0)\in (H\sp R)\sperp$ and $(0, \eta)\in (H\sp S)\sperp$,
Claim~\ref{claim:1} is proved.
\end{proof}
Because  $H\sperp/H$ is $3$-elementary,
we have $(0, \gamma)\notin H\sperp$
by~\eqref{eq:HcapS}.
Hence we obtain
\begin{equation}\label{eq:ne}
 (H\sp S)\sperp=(H\sperp)\sp S \ne \Disc{S}.
\end{equation}
By~\eqref{eq:inKer}, $H\sp S$ is either $0$ or $\Ker m\sb 3$.
If $H\sp S=0$, then $(H\sp S)\sperp=\Disc{S}$ and we get a contradiction to~\eqref{eq:ne}.
Suppose that $H\sp S=\Ker m\sb 3$.
Then $(H\sp S)\sperp\supset \Im m\sb 3$, and hence 
$(H\sp S)\sperp= \Im m\sb 3$ by~\eqref{eq:ne}.
In particular,
we have
\begin{equation}\label{eq:one}
\log \sb 3 |(H\sp S)\sperp|=1.
\end{equation}
Since $|(H\sp R)\sperp|=3^{10}/ |H\sp R|$ and $H\cong H\sp R$  by~\eqref{eq:HcapS},
we see that 
$$
2\sigma = \log \sb 3 |H\sperp/H|=\log \sb 3 |(H\sp R)\sperp| + \log \sb 3 |(H\sp S)\sperp| -\log \sb 3 |H| 
=10-2 \log \sb 3 |H|  +1
$$
is odd by~\eqref{eq:one}, which is absurd.
Therefore $\Disc{S}\not\cong \Z/9\Z$.

Because $S$ is an even lattice, 
the classification of indefinite lattices of rank $2$~(\cite[Chapter 15, Section 3]{CS})
implies the following:
\begin{eqnarray*}
\Disc{S}=0 & \;\;\Longrightarrow\;\; & S\cong U(1), \\
\Disc{S}\cong \F\sb 3\sp 2  & \;\;\Longrightarrow\;\; & S\cong U(3), \\
\Disc{S}\cong (\Z/9\Z)\sp 2 & \;\;\Longrightarrow\;\; & S\cong U(9). 
\end{eqnarray*}

Next we assume  $S\cong U(9)$, and derive a contradiction.
Note that $\Ker m\sb 3$ is 
generated by 
$$
3\bar e\dual=f/3 \;\bmod S \quad\textrm{and}\quad 3\bar f\dual=e/3 \;\bmod S.
$$
By~\eqref{eq:inKer2}, we have
\begin{equation}\label{eq:inF}
H\sperp\;\;\subset\;\; \Disc{R} \oplus \Ker m\sb 3.
\end{equation}
Then $H$ is also contained in $ \Disc{R} \oplus \Ker m\sb 3$.
Suppose that $H$ is generated by
$$
g\sp{(\nu)}=\xi\sp{(\nu)}\sb 1 \gamma\sb 1 +\cdots +  \xi\sp{(\nu)}\sb{10} \gamma\sb{10} + 
\eta\sp{(\nu)}\sb 1 (3\bar e\dual ) + \eta\sp{(\nu)}\sb 2 (3\bar f\dual)
\qquad(\nu=1, \dots, r)
$$
where $\xi\sp{(\nu)}\sb i, \eta\sp{(\nu)}\sb j\in \F\sb 3$.
We put
$$
M:=\left[
\begin{array}{cccc|cc}
\xi\sb{1}\sp{(1)} & \dots & \dots & \xi\sb{10}\sp{(1)} &\eta\sp{(1)}\sb 2 & \eta\sp{(1)}\sb 1\\
  & \dots & \dots &  & & \\
  & \dots & \dots &  & & \\
\xi\sb{1}\sp{(r)} & \dots & \dots & \xi\sb{10}\sp{(r)} &\eta\sp{(r)}\sb 2 & \eta\sp{(r)}\sb 1\\
\end{array}
\right].
$$
From~\eqref{eq:discform10At} and~\eqref{eq:discformUm},  an element
$$
x\sb 1\gamma\sb 1 +\cdots + x\sb{10} \gamma\sb{10} + 
y\sb 1 \bar e\dual  +  y\sb 2 \bar f\dual \qquad
(x\sb 1, \dots, x\sb{10} \in \F\sb 3, \;\; y\sb 1, y\sb 2 \in \Z/9\Z)
$$
of $\Disc{R} \oplus \Disc{S}$ is contained in $H\sperp$ if and only if 
the vector
$\vx:=[x\sb 1, \dots, x\sb{10}, y\sb 1, y\sb 2]$ satisfies the  equation
\begin{equation}\label{eq:lineq}
M\cdot {}\sp T \vx \equiv \vz \;\bmod  3.
\end{equation}
We consider~\eqref{eq:lineq} as a system of linear equations over $\F\sb 3$.
The property~\eqref{eq:inF} of $H\sperp$ implies that every solution of~\eqref{eq:lineq} in $\F\sb 3$
must satisfy
\begin{equation}\label{eq:yzero}
y\sb 1 =y\sb 2=0.
\end{equation}
Because of~\eqref{eq:HcapS} and hence $H \cong H^R$,
we can choose generators $g\sp{(1)}$, \dots, $g\sp{(r)}$ of $H$ in such a way that,
after suitable permutations of $10$ coordinates of $\Disc{R}=\F\sb 3\sp{10}$ if necessary,
the $r\times 12$ matrix $M$ is of the form
$$
M=\left[
\begin{array}{ccc|c}
 &&& \\
 &I\sb r &&* \\
 &&& \\
\end{array}
\right], 
$$
where $I\sb r$ ($r\le 10$) is a diagonal matrix
whose diagonal entries are $1$.
Now non-zero elements of 
the subgroup of $H\cong H^R $ of $ H^{\perp}$
should be solutions of~\eqref{eq:lineq} in $\F\sb 3$, but
do not satisfy~\eqref{eq:yzero}.
Thus we get a contradiction. 

Hence $S$ is isomorphic to $U(1)$ or $U(3)$.
If $S\cong U(1)$, then 
$2\sigma =10 -2 \dim\sb{\F\sb 3} H\le 10 $,
while 
if $S\cong U(3)$, then 
$2\sigma =10+2  -2 \dim\sb{\F\sb 3} H\le 12 $.
\end{proof}
\begin{example}\label{example:quartic}
Let $[w:x:y:z]$ be homogeneous coordinates of $\P\sp 3$.
For homogeneous polynomials $f(y, z)$, $g(y, z)$ and $h(y, z)$ of degrees $3$, $3$ and $4$,
we consider the quartic surface $X$ defined in $\P\sp 3$ by
$$
 w^3 y + x^3z + w f(y, z) + x g(y, z) + h(y, z) = 0.
$$
When $X$ is smooth,
$X$ is a supersingular $K3$ surface,
because $X$ contains a configuration of lines 
as in ~\cite[Section 6]{Shimada92}.
%
%
It was shown in~\cite[Section 6]{Shimada92} that, 
when $f$, $g$ and $h$ are general,
the Artin invariant of $X$ is $6$,
and hence the orthogonal complement $R\sperp$ of the sublattice $R\st \NS (X)$
generated by the classes of the lines $C\sb 1, D\sb 1, \dots, C\sb{10}, D\sb{10}$
(10 pairs of intersecting lines in 10 singular fibres of type IV)
is isomorphic to $U(3)$.
When $f$, $g$  are general and $h=0$,
the Artin invariant of $X$ is $5$ by~\cite[Section 4]{Shioda77}.
In this case,
the line $\ell$  defined by $w=x=0$ is contained in $X$.
Since $\ell^2=-2$ and $[\ell]\in R\sperp$,
$R\sperp$ is isomorphic to $U(1)$.
%
%
\end{example}
\section{Proof of Theorem~\ref{thm:U3}}\label{sec:U3}
The discriminant group $D$ of $\Z[10\At]\oplus U(3)$ is equal to
$$
\F\sb 3 \gamma\sb 1\oplus \cdots \oplus \F\sb 3 \gamma\sb{10} \oplus \F\sb 3 \bar e\dual  \oplus \F\sb 3 \bar f\dual,
$$
and the discriminant form $q$ of $\Z[10\At]\oplus U(3)$ is given by
$$
q(x\sb 1, \dots, x\sb{10}, y\sb 1, y\sb 2)=-2(x\sb 1\sp 2 + \cdots + x\sb{10}\sp 2 )/3 + 2 y\sb 1 y\sb 2 /3 
\;\;\in\;\;\Q/2\Z.
$$
We consider subgroups of $D$ as ternary codes.
Recall from \S\ref{sec:disc} that
 the Hamming weight of a word $\vx=(x\sb 1, \dots, x\sb{10})\in \Disc{\Z[10\At]}$ is defined by 
$$
\wt (\vx):=|\set{i}{x\sb i\ne 0}|.
$$
Then a ternary code $\CCC\st D$ is isotropic with respect to 
$q$ if and only if
\begin{equation}\label{eq:isot3}
\wt (\vx)\equiv y\sb 1 y\sb 2 \;\bmod 3\;\;\textrm{for any \;\;$(\vx, y\sb 1, y\sb 2)\in \CCC$}
\end{equation}
holds.
A ternary code $\CCC\st D$ satisfying~\eqref{eq:isot3} is therefore called an \emph{isotropic code}.
For an isotropic code $\CCC$,
we denote by $N\sb{\CCC}$ the overlattice of $\Z[10\At]\oplus U(3)$
corresponding to $\CCC$ by Proposition~\ref{prop:nikulin}.
By Theorem~\ref{thm:NRS},
$N\sb{\CCC}$ is isomorphic to the lattice $N\sb{3, \sigma}$, where $\sigma =6-\dim \CCC$.

It is easy to see that the following conditions for an isotropic code $\CCC$ are equivalent:
\begin{itemize}
\item[(i)] $\wt (\vx)>0 \;\;\; \textrm{for any non-zero word \;$(\vx, y\sb 1, y\sb 2)\in \CCC$}$,
\item[(ii)] $U(3)$ is primitive in $N\sb{\CCC}$, and 
\item[(iii)] $\Z [10\At]\sperp=U(3)$ in  $N\sb{\CCC}$.
\end{itemize}
We say that an isotropic code $\CCC$ is \emph{admissible} if $\CCC$ satisfies the conditions above.
Let $h=a e+ bf$ be a vector of $U(3)$ with $a\ge 1$ and  $b\ge 1$.
We have
$h^2=6ab$.
\begin{lemma}\label{lemma:abU3}
Let $\CCC$ be an admissible isotropic code.

{\rm (1)} There exists a vector $u\in N\sb{\CCC}$ satisfying $hu=1$ or $2$ and $u^2=0$ if and only if
the following hold:
\begin{itemize}
\item[{\rm ($\alpha$)}] $a=b=1$,  and
\item[{\rm ($\beta$)}] there exists  $(\vx, y\sb 1, y\sb 2)\in \CCC$ such that $\wt (\vx)=1$.
\end{itemize}

{\rm (2)} The set of roots
$\Roots (h\sperp):=\shortset{r\in N\sb{\CCC}}{rh=0,\; r^2=-2}$ in $h\sperp$
is strictly larger than $\Roots(\Z[10\At])=\{ \pm c\sb i, \pm d\sb i, \pm (c\sb i + d\sb i)\}$
if and only if one of the following holds:
\begin{itemize}
\item[{\rm (a)}]
there exists  $(\vx, y\sb 1, y\sb 2)\in \CCC$ such that $\wt (\vx)=3$ and $y\sb 1=y\sb 2=0$, or 
\item[{\rm (b)}]
$a=b$,
 and 
there exists  $(\vx, y\sb 1, y\sb 2)\in \CCC$ such that $\wt (\vx)=2$, or
\item[{\rm (c)}]
{\rm (}$a=2b\;\;\textrm{or}\;\;b=2a${\rm )}
 and 
there exists  $(\vx, y\sb 1, y\sb 2)\in \CCC$ such that $\wt (\vx)=1$.
\end{itemize}
\end{lemma}
\begin{proof}
We prove (1) first. Suppose that a vector
\begin{equation}\label{eq:u}
u=r\sb u  + \eta\sb 1 e\dual + \eta\sb 2 f\dual \quad (r\sb u \in \Z[10\At]\dual, \;\; \eta\sb 1, \eta\sb 2\in \Z)
\end{equation}
of $N\sb{\CCC}$ satisfies $hu=1$ or $2$ and $u^2=0$.
Then we have
\begin{eqnarray}
&& a \eta\sb 1 + b \eta\sb 2 = 1 \;\textrm{or}\; 2,  \label{eq:1or2} \\
&& r\sb u\sp 2 + 2 \eta\sb 1 \eta\sb 2/3 =0. \label{eq:rueta}
\end{eqnarray}
Note that $(\eta\sb 1, \eta\sb 2)\not\equiv (0,0) \bmod 3$ by~\eqref{eq:1or2}.
Since $\CCC$ is admissible, we have $r\sb u\ne 0$ by , and hence
$\eta\sb 1\eta\sb 2>0$ by~\eqref{eq:rueta}.
From~\eqref{eq:1or2},
we obtain
$$
a=b=1,  \qquad \eta\sb 1=\eta\sb 2=1,
$$
and hence, from~\eqref{eq:rueta}, we have
$$
r\sb u\sp 2 =-2/3.
$$
By~\eqref{eq:wtr}, the word 
$$
\bar u=u \; \bmod  (\Z[10\At]\oplus U(3))=(\bar{r}\sb u, \bar\eta\sb 1, \bar\eta\sb 2)\qquad
(\textrm{where $\bar{r}\sb u= r\sb u \bmod \Z[10\At]$})
$$
 of $\CCC$ has the  property $\wt (\bar{r}\sb u)=1$.

\par
Conversely, suppose that $a=b=1$ and that 
there exists a word  $(\bar r , y\sb 1, y\sb 2)\in \CCC$ such that $\wt (\bar r)=1$.
Replacing $(\bar r , y\sb 1, y\sb 2)$ by $(-\bar r , -y\sb 1, -y\sb 2)$
if necessary, we can assume that $y\sb 1=y\sb 2=1$ by
~\eqref{eq:isot3}.
Then, by
~\eqref{eq:wtrinv}, there exists a vector
$$
u=r+ e\dual + f\dual\qquad(r\in \Z[10\At]\dual)
$$
in $N\sb{\CCC}$ satisfying  $r^2=-2/3$.
This vector $u$ satisfies $hu=2$ and $u^2=0$.
Thus the assertion (1) is proved.

\par
We now prove (2).
Suppose that a vector $u\in N\sb{\CCC}$ given by~\eqref{eq:u} satisfies $hu=0$, $u^2=-2$ and $u\notin \Roots(\Z[10\At])$.
Then we have
\begin{eqnarray}
&& a \eta\sb 1 + b \eta\sb 2 = 0,  \label{eq:eta0} \\
&& r\sb u\sp 2 + 2 \eta\sb 1 \eta\sb 2/3 =-2. \label{eq:rueta2}
\end{eqnarray}
Suppose that  $\eta\sb 1=0$ or $\eta\sb 2 =0$.
Then~\eqref{eq:eta0} implies $\eta\sb 1=\eta\sb 2=0$ and hence 
$\wt (\bar{r}\sb u)\equiv 0  \bmod  3$ holds because $\CCC$ is isotropic.
By~\eqref{eq:wtr} and~\eqref{eq:rueta2}, we have $\wt (\bar{r}\sb u)\le 3$.
If $\wt (\bar{r}\sb u)=0$, then $u=r\sb u$ is contained in $\Roots (\Z[10\At])$.
Hence we have $\wt (\bar{r}\sb u)=3$,
and therefore the condition (a) is satisfied.
Suppose that $\eta\sb 1 \ne 0$ and $\eta\sb 2\ne 0$.
By~\eqref{eq:eta0}, we have $\eta\sb 1\eta\sb 2<0$.
By~\eqref{eq:wtr} and~\eqref{eq:rueta2},
we see that the pair $(\eta\sb 1\eta\sb 2, \wt (\bar{r}\sb u))$ is either $(-1, 2)$ or $(-2, 1)$.
In the former case, we have $a=b$ by~\eqref{eq:eta0} and hence (b) is satisfied.
In the latter case, we have $a=2b$ or $b=2a$ by~\eqref{eq:eta0} and hence (c) is satisfied.

\par
Conversely, suppose that (a) is fulfilled.
Using~\eqref{eq:wtrinv}, we have a lift 
$$
u=r+0+0\in N\sb{\CCC}\qquad (r\in \Z [10\At]\dual)
$$
of the word $(\bar r , 0, 0)\in \CCC$ with $\wt(\bar r)=3$
such  that $r^2=-2$. Then $u\in \Roots(h\sperp)\sm \Roots (\Z[10\At])$.
Suppose that (b) is satisfied.
A vector 
$$
u=r+e\dual-f\dual\;\;\in\;\; N\sb{\CCC}
$$
with $\wt(\bar r)=2$ and $r^2=-4/3$
satisfies  $u\in \Roots(h\sperp)\sm \Roots (\Z[10\At])$.
Suppose that (c) is satisfied and assume that $a=2b$.
A vector 
$$
u=r+e\dual-2 f\dual\;\;\in\;\; N\sb{\CCC}
$$
with $\wt(\bar r)=1$ and $r^2=-2/3$
satisfies   $u\in \Roots(h\sperp)\sm \Roots (\Z[10\At])$.
Thus the assertion (2) is proved.
\end{proof}
\begin{proof}[Proof of Theorem~\ref{thm:U3}]
The implication $\textrm{(iii)}\Longrightarrow\textrm{(ii)}$ is obvious.
Since every vector $h$ of $U(3)$ satisfies $h^2\equiv 0 \bmod 6$,
the implication $\textrm{(ii)}\Longrightarrow\textrm{(i)}$ is also obvious.
Using computer,
we can prove the following Claim~\ref{claim:sevencodes}.
See Remark~\ref{remark:sevencodes} and Table~\ref{table:sevencodes7}.
\begin{claim}\label{claim:sevencodes}
There exists an isotropic admissible code $\CCC\st D$ of dimension $5$ with the following property:
\begin{equation}\label{eq:theproperty}
\parbox{10cm}{
every non-zero word $(\vx, y\sb 1, y\sb 2)\in \CCC$ satisfies the following:
(i) $\wt (\vx)\ge 3$, and (ii) if $\wt (\vx)=3$, then $(y\sb 1, y\sb 2)\ne(0, 0)$.
}
\end{equation}
\end{claim}
We now prove $\textrm{(i)}\Longrightarrow\textrm{(iii)}$
Suppose that an integer $d=6m \;\; (m\in \Z\sb{>0})$ is given.
Let $X$ be a supersingular $K3$ surface in characteristic $3$
with Artin invariant $\sigma\le 6$.
For the basis $e, f$ of $U(3)$ at the end of Section 2, we put
$$
h:=e+mf.
$$
Then $h^2=d$.
Let $\CCC (\sigma)$ be a linear subspace of the code $\CCC$ 
in Claim~\ref{claim:sevencodes}
with $\dim \CCC(\sigma)=6-\sigma$.
Since $\CCC(\sigma)$ is isotropic,
the corresponding overlattice $N\sb{\CCC(\sigma)}$
of $\Z[10\At]\oplus U(3)$ is isomorphic to $N\sb{3, \sigma}$ by Theorem~\ref{thm:NRS}.
Hence there exists an isometry
$$
\phi : N\sb{\CCC(\sigma)}\isom \NS(X)
$$
by  Theorem~\ref{thm:ARS}.
Since every word of $\CCC (\sigma)$ satisfies the conditions (i) and (ii) in~\eqref{eq:theproperty},
Lemma~\ref{lemma:abU3} implies that
there exist no vectors $u$ in $N\sb{\CCC(\sigma)}$ satisfying $hu=1$ or $2$ and $u^2=0$,
and that the set of roots in the orthogonal complement $h\sperp$ of $h$ in $N\sb{\CCC(\sigma)}$
coincides with $\Roots(\Z[10\At])$.
By Proposition~\ref{prop:h} and Remark~\ref{rem:degree8}, 
we can choose the isometry 
$\phi : N\sb{\CCC(\sigma)}\isom \NS(X)$ in such a way that $\phi (h)$ is the class $[L]$
of a line bundle $L$ very ample modulo $(-2)$-curves
such that $\Phi\sb{|L|}$ induces a contraction $\rho\sb L : X\to Y\sb{\polX}$ 
of an $ADE$-configuration of $(-2)$-curves of type $10\At$.
Since $\CCC (\sigma)$ is admissible,
we see that $R\sperp \sb{\polX}\subset \NS (X)$ is isomorphic to $U(3)$.
Thus $X$ admits a polarization $L$ of degree $d$ with the hoped-for properties.
\end{proof}
\begin{table}
$$
\CCC\sb 1\quad : \quad 
\left [\begin {array}{cccccccccccc} 
1&0&0&0&0&0&0&0&1&1&0&1\\
0&1&0&0&0&0&0&1&0&1&2&0\\
0&0&1&0&0&0&1&0&1&0&2&0\\
0&0&0&1&0&0&1&1&0&0&0&1\\
0&0&0&0&1&0&1&1&1&1&1&2
\end {array}\right ]
$$
\par\smallskip
$$
\CCC\sb 2\quad : \quad \left [\begin {array}{cccccccccccc} 
1&0&0&0&0&0&0&0&1&1&0&1\\
0&1&0&0&0&0&0&1&0&1&2&0\\
0&0&1&0&0&0&1&0&1&0&2&0\\
0&0&0&1&0&0&1&1&0&0&0&1\\
0&0&0&0&1&1&1&2&2&1&0&0
\end {array}\right ]
$$
\par\smallskip
$$
\CCC\sb 3\quad : \quad \left [\begin {array}{cccccccccccc} 
1&0&0&0&0&0&0&0&1&1&0&1\\
0&1&0&0&0&0&1&1&0&0&0&1\\
0&0&1&0&0&1&0&1&0&1&2&2\\
0&0&0&1&0&1&1&0&1&0&2&2\\
0&0&0&0&1&1&1&2&2&1&0&1
\end {array}\right ]
$$
\par\smallskip
$$
\CCC\sb 4\quad : \quad \left [\begin {array}{cccccccccccc} 
1&0&0&0&0&0&0&0&1&1&0&1\\
0&1&0&0&0&0&1&1&0&0&0&1\\
0&0&1&0&0&1&0&1&0&1&2&2\\
0&0&0&1&0&1&1&0&1&0&2&2\\
0&0&0&0&1&1&2&2&2&2&2&0
\end {array}\right ]
$$
\par\smallskip
$$
\CCC\sb 5\quad : \quad \left [\begin {array}{cccccccccccc} 
1&0&0&0&0&0&0&1&1&1&1&1\\
0&1&0&0&0&0&1&0&1&1&2&2\\
0&0&1&0&0&1&0&1&0&1&2&2\\
0&0&0&1&0&1&1&0&0&1&1&1\\
0&0&0&0&1&1&1&1&1&1&0&0
\end {array}\right ]
$$
\par\smallskip
$$
\CCC\sb 6\quad : \quad \left [\begin {array}{cccccccccccc} 
1&0&0&0&0&0&0&1&1&1&1&1\\
0&1&0&0&0&0&1&0&1&1&2&2\\
0&0&1&0&0&1&0&1&0&1&2&2\\
0&0&0&1&0&1&1&0&0&1&1&1\\
0&0&0&0&1&1&2&2&1&0&1&2
\end {array}\right ]
$$
\par\smallskip
$$
\CCC\sb 7\quad : \quad \left [\begin {array}{cccccccccccc} 
1&0&0&0&0&0&1&1&1&1&1&2\\
0&1&0&0&0&1&0&1&1&2&2&1\\
0&0&1&0&0&1&1&0&2&1&2&1\\
0&0&0&1&0&1&1&2&0&2&1&2\\
0&0&0&0&1&1&2&1&2&0&1&2
\end {array}\right ]
$$
\par\medskip
\caption{Bases of the codes 
$\CCC\sb 1$, \dots, $\CCC\sb 7$}\label{table:sevencodes7}
\end{table}
\begin{table}
\begin{eqnarray*}
\we(\CCC\sb 1) &=& 1+12\,{z}^{3}+18\,{z}^{4}+36\,{z}^{5}+108\,{z}^{6}+36\,{z}^{7}+18\,{z}^{8}+14\,{z}^{9}, \\
\we(\CCC\sb 2) &=& 1+8\,{z}^{3}+10\,{z}^{4}+24\,{z}^{5}+86\,{z}^{6}+40\,{z}^{7}+30\,{z}^{8}+40\,{z}^{9}+4\,{z}^{10},\\
\we(\CCC\sb 3) &=& 1+4\,{z}^{3}+8\,{z}^{4}+24\,{z}^{5}+94\,{z}^{6}+44\,{z}^{7}+30\,{z}^{8}+36\,{z}^{9}+2\,{z}^{10}, \\
\we(\CCC\sb 4) &=& 1+6\,{z}^{3}+6\,{z}^{4}+18\,{z}^{5}+102\,{z}^{6}+42\,{z}^{7}+36\,{z}^{8}+26\,{z}^{9}+6\,{z}^{10}, \\
\we(\CCC\sb 5) &=& 1+30\,{z}^{4}+60\,{z}^{6}+120\,{z}^{7}+20\,{z}^{9}+12\,{z}^{10}, \\
\we(\CCC\sb 6) &=& 1+18\,{z}^{4}+18\,{z}^{5}+96\,{z}^{6}+36\,{z}^{7}+36\,{z}^{8}+38\,{z}^{9}, \\
\we(\CCC\sb 7) &=& 1+72\,{z}^{5}+60\,{z}^{6}+90\,{z}^{8}+20\,{z}^{9}.
\end{eqnarray*}
\par\medskip
\caption{Weight-enumerators}\label{table:wes}
\end{table}
\begin{remark}\label{remark:sevencodes}
Let $G$ denote the group of linear automorphisms of $D\cong \F\sb 3 \sp{10}\oplus \F\sb 3 \sp 2 $ 
generated by
\begin{eqnarray*}
(x\sb 1, \dots, x\sb{10}, y\sb 1, y\sb 2 ) &\mapsto & 
(x\sb{\sigma(1)}, \dots, x\sb{{\sigma(10)}}, y\sb {\tau(1)}, y\sb {\tau(2)} )
\qquad (\sigma \in \SSSS\sb{10}, \tau\in \SSSS\sb 2), \quad\textrm{and}\\
(x\sb 1, \dots, x\sb{10}, y\sb 1, y\sb 2 ) &\mapsto & 
((-1)\sp{\alpha\sb 1}x\sb 1, \dots, (-1)\sp{\alpha\sb {10}}x\sb {10}, (-1)\sp\beta y\sb 1, (-1)\sp\beta y\sb 2 )\\
&& \qquad (\alpha\sb 1, \dots, \alpha\sb{10}\in \F\sb 2, \beta\in \F\sb 2 ).
\end{eqnarray*}
Note that, if $\CCC\st D$ is an isotropic admissible code,
then so is $g (\CCC)$ for any $g\in G$.
We  define the weight enumerator of a ternary code $\CCC$ by
$$
\we (\CCC):=\sum\sb{(\vx, y\sb 1, y\sb 2)\in \CCC} z\sp{\wt (\vx)}.
$$
Using computer, we have proved that there exist at least  seven isomorphism  classes
of isotropic admissible codes of dimension $5$ with the property~\eqref{eq:theproperty}.
The representative codes $\CCC\sb 1, \dots, \CCC\sb 7$ of these classes are given in
Table~\ref{table:sevencodes7}.
Their weight-enumerators are given in Table~\ref{table:wes}.
\end{remark}
\begin{corollary}\label{cor:seven}
Let $X$ be a supersingular $K3$ surface in characteristic $3$
with Artin invariant $1$.
Then there exist at least seven  line bundles $L\sb 1, \dots, L\sb 7$ of degree $6$
on $X$ that are mutually non-isomorphic and that induce  contractions
of $10 A\sb 2$-configurations of $(-2)$-curves on $X$.
\end{corollary}
See Example~\ref{example:sigmaG0}.
\section{Proof of Theorem~\ref{thm:U1}}\label{sec:U1}
The proof of Theorem~\ref{thm:U1} is similar to and simpler than
that of Theorem~\ref{thm:U3}.
\par
\medskip
The discriminant group $D$ of $\Z[10\At]\oplus U(1)$ is equal to
$$
\F\sb 3 \gamma\sb 1\oplus \cdots \oplus \F\sb 3 \gamma\sb {10}.
$$
A ternary code $\CCC\st D$ is isotropic with respect to 
the discriminant form $q$ of $\Z[10\At]\oplus U(1)$ if and only if
\begin{equation}\label{eq:isot1}
\wt (\vx)\equiv 0 \;\bmod 3\;\; \textrm{for any \;\;$\vx\in \CCC$}
\end{equation}
holds.
For an isotropic code $\CCC$,
we denote by $N\sb{\CCC}$ the overlattice of $\Z[10\At]\oplus U(1)$
corresponding to $\CCC$.
By Theorem~\ref{thm:NRS},
$N\sb{\CCC}$ is isomorphic to the lattice $N\sb{3, \sigma}$, where $\sigma =5-\dim \CCC$.

Let $h=a e+ bf$ be a vector of $U(1)$ with $a\ge 1$ and  $b\ge 1$.
We have $h^2=2ab$.
\begin{lemma}\label{lemma:ab}
Let $\CCC$ be an  isotropic code in $D\cong \F\sb 3 \sp{10}$.

{\rm (1)} There exists a vector $u\in N\sb{\CCC}$ satisfying $hu=1$ or $2$ and 
$u^2=0$ if and only if $a\le 2$ or $b\le 2$.

{\rm (2)} The set of roots
$\Roots (h\sperp):=\shortset{r\in N\sb{\CCC}}{rh=0,\; r^2=-2}$ in $h\sperp$ 
is strictly larger than $\Roots(\Z[10\At])$
if and only if one of the following holds;
\begin{itemize}
\item[{\rm (a)}]
there exists  $\vx\in \CCC$ such that $\wt (\vx)=3$, or 
\item[{\rm (b)}]
$a=b$.
\end{itemize}
\end{lemma}
\begin{proof}
We prove (1) first. Suppose that a vector
\begin{equation}\label{eq:uU1}
u=r\sb u  + \eta\sb 1 f + \eta\sb 2 e \quad (r\sb u \in \Z[10\At]\dual, \;\; \eta\sb 1, \eta\sb 2\in \Z)
\end{equation}
of $\Z[10\At]\dual \oplus U(1)\dual=\Z[10\At]\dual \oplus U(1)$ satisfies $hu=1$ or $2$ and $u^2=0$.
Then we have
\begin{eqnarray}
&& a \eta\sb 1 + b \eta\sb 2 = 1 \;\textrm{or}\; 2,  \label{eq:1or2U1} \\
&& r\sb u\sp 2 + 2 \eta\sb 1 \eta\sb 2 =0. \label{eq:ruetaU1}
\end{eqnarray}
By~\eqref{eq:ruetaU1},
we have $\eta\sb 1\eta\sb 2 \ge 0$.
Using~\eqref{eq:1or2U1}, we have $a\le 2$ or $b\le 2$.
Conversely, if $a\le 2$, then $u=f$ satisfies  $hu=a=1$ or $2$ and $u^2=0$.
Thus (1) is proved.

\par
Next we prove (2).
Suppose that a vector $u$ given by~\eqref{eq:uU1} satisfies $hu=0$, $u^2=-2$ and $u\notin \Roots(\Z[10\At])$.
Then we have
\begin{eqnarray}
&& a \eta\sb 1 + b \eta\sb 2 = 0,  \label{eq:eta0U1} \\
&& r\sb u\sp 2 + 2 \eta\sb 1 \eta\sb 2 =-2. \label{eq:rueta2U1}
\end{eqnarray}
Because $\CCC$ is isotropic, 
$\wt (\bar{r}\sb u)\equiv 0  \bmod  3$ holds.
If $\eta\sb 1=0$ or $\eta\sb 2 =0$, then~\eqref{eq:eta0U1} implies $\eta\sb 1=\eta\sb 2=0$.
By~\eqref{eq:wtr} and~\eqref{eq:rueta2U1}, we have $\wt (\bar{r}\sb u)\le 3$.
If $\wt (\bar{r}\sb u)=0$, then $u=r\sb u$ is contained in $\Roots (\Z[10\At])$.
Hence we have $\wt (\bar{r}\sb u)=3$,
and therefore the condition (a) is satisfied.
Suppose that $\eta\sb 1 \ne 0$ and $\eta\sb 2\ne 0$.
By~\eqref{eq:eta0U1}, we have $\eta\sb 1\eta\sb 2<0$.
By~\eqref{eq:rueta2U1},
we have  $r\sb u=0$ and $\eta\sb 1\eta\sb 2=-1$,
and hence $a=b$ follows from~\eqref{eq:eta0U1}.

\par
Conversely, suppose that (a) is fulfilled.
Using~\eqref{eq:wtrinv}, we have a lift 
$$
u=r+0+0\in N\sb{\CCC}\qquad (r\in \Z [10\At]\dual)
$$
of the word $\bar r \in \CCC$ with $\wt(\bar r)=3$
such  that $r^2=-2$. Then $u$ is contained in $\Roots(h\sperp)\sm \Roots (\Z[10\At])$.
Suppose that (b) is satisfied.
The vector 
$$
u=e-f\;\;\in\;\; N\sb{\CCC}
$$
satisfies   $u\in \Roots(h\sperp)\sm \Roots (\Z[10\At])$.
\end{proof}
In order to prove Theorem~\ref{thm:U1}, it is therefore enough to show the following:
\begin{claim}
There exists an isotropic code $\CCC\st D\cong \F\sb 3\sp{10}$
of dimension $4$
such that $\wt (\vx)\ge 6$ holds for any $\vx \in \CCC$.
\end{claim}
The code  $ \CCC$ generated by the row vectors  of 
$$
\left [\begin {array}{cccccccccc} 
1&0&0&0&0&1&1&1&1&1\\
0&1&0&0&1&0&1&1&2&2\\
0&0&1&0&1&1&0&2&1&2\\
0&0&0&1&1&1&2&0&2&1
\end {array}\right ]
$$
satisfies  $\wt (\vx)\ge 6$  for any $\vx \in \CCC$.
The weight-enumerator $\sum\sb{\vx\in \CCC} z\sp{\wt (\vx)}$
of this  code $\CCC$  is 
$$
1+60 z^6+20 z^9.
$$
\begin{remark}
The code $\CCC$ above  is obtained as 
a subcode of the extended ternary Golay code in $\F\sb 3\sp{12}$.  See~\cite[Chapter 5, Section 2]{E}.
\end{remark}
\section{Proof of Theorem~\ref{thm:insep}}\label{sec:insep}
Let $(X, L)$ be a polarized $K3$ surface of degree $6$.
Then $Y\sb{(X, L)}$ is  
a complete intersection of multi-degree $(2, 3)$ in $\Pf$ by \cite[Theorem 6.1]{SD}.
Let $\tQ\sb{(X, L)}$ denote the unique quadric hypersurface in $\Pf$
containing $Y\sb{(X, L)}$.
\begin{proposition}
Suppose that 
$\RRR\sb{(X, L)}=10\At$ and
$R\sb{\polX}\sperp \cong U(3)$.
Then $\tQ\sb{(X, L)}$ is a cone over a non-singular quadric surface $Q=\Pll$,
and $Y\sb{(X, L)}$ does not pass through the vertex $P$ of the cone $\tQ\sb{(X, L)}$.
\end{proposition}
\begin{proof}
By the assumption,
$R\sb{\polX}\sperp$ is generated by the numerical equivalence classes  $[E]$ and $[F]$ 
of divisors $E$ and $F$ satisfying
\begin{equation}\label{eq:EF}
E^2=F^2=0, \quad EF=3, \quad [L]=[E]+[F].
\end{equation}
By the Riemann-Roch theorem, we can assume that $E$ and $F$ are effective.
Suppose that $|E|$ has a fixed component.
Let $M+\Gamma$ be a general member of $|E|$,
where $\Gamma$ is the fixed part of $|E|$.
Because $\rho\sb L: X\to Y\sb{\polX}$ is birational,
$\rho\sb L$ induces a birational map from $M$ to $\rho\sb L (M)$.
Note that $\rho\sb L (M+\Gamma)$ is a cubic curve.
If $\rho\sb L (\Gamma)$ is of dimension $1$,
then $\rho\sb L (M)$ is a line or a plane conic,
and hence contradicts $\dim |M|>0$.
Therefore, $\rho\sb L$ contracts every irreducible component of
$\Gamma$ to a point,
and hence $[\Gamma]\in R\sb{\polX}$.
From $[E]\in R\sb{\polX}\sperp$,
we obtain $E\Gamma=0$ and hence 
$M^2=E^2+\Gamma^2<0$.
Thus we get a contradiction again.
Hence  $|E|$ has no fixed components.
In particular, $E$ is nef.
Since  $\rho\sb L$ is birational and $E$ is primitive in $R\sb{\polX}\sperp$ 
(being part of its basis), 
a general member $E$ of $|E|$ is mapped by $\rho\sb L$ birationally 
to a plane cubic curve in $\Pf$.
Therefore  a general member of $|E|$ is irreducible,
and hence $|E|$ is a (quasi-)elliptic pencil by~\cite[Proposition 0.1]{Nikulin79}.
Therefore the quadric hypersurface $\tQ\sb{(X, L)}$  contains a one-dimensional family 
$\{\Pi\sp E\sb t\}$
of  planes
such that 
$$
|E|=\{\rho\sb L \sp *  (\Pi\sp E\sb t\cap Y\sb{\polX})\}.
$$
Hence $\tQ\sb{(X, L)}$ is singular.
Since $\tQ\sb{(X, L)}$ contains two irreducible families $\{\Pi\sp E\sb t\}$ and $\{\Pi\sp F\sb t\}$ of planes 
corresponding to $|E|$ and $|F|$,
we have $\dim \Sing \tQ\sb{(X, L)}=0$,
and $\tQ\sb{(X, L)}$ is a cone over a non-singular quadric surface $Q=\Pll$.
If $Y\sb{(X, L)}$  passed through the vertex $P$ of the cone $\tQ\sb{(X, L)}$,
then the linear system $|E|$
would have a fixed component that is contracted to the point $P$.
Hence  $P$ is not contained in $Y\sb{(X, L)}$. 
\end{proof}
Note that a non-ordered pair of the numerical equivalence classes $[E]$ and $[F]$  in $R\sb{\polX}\sperp$ 
satisfying~\eqref{eq:EF} is unique.
The following has been  shown in the proof above:
\begin{corollary}\label{cor:EFnef}
The divisors $E$ and $F$ are nef.
The complete linear systems $|E|$ and $|F|$ are {\rm (}quasi-{\rm )}elliptic pencils.
\end{corollary}
We denote by
$$
\pi\sb P : Y\sb{\polX}\to Q=\Pll
$$
the projection from the vertex $P$ of the cone $\tQ\sb{(X, L)}$.
Let $x$ and $y$ be affine coordinates of the two factors of $\Pll$.
The surface $Y\sb{(X, L)}$ is defined by an equation 
\begin{equation}\label{eq:F}
\Psi\;:=\; W^3\,+\, a(x, y)\, W^2\,+ \,b(x, y)\, W \, + \,c(x, y)\;=\;0,
\end{equation} 
where $W$ is a fiber coordinate of the affine line bundle $\tQ\sb{(X, L)}\sm \{P\}\cong L\sb Q(1,1)$
on $Q=\Pll$,
and $a$, $b$, $c$ are polynomials of degrees $1$, $2$ and  $3$, respectively.
\par
\medskip
Let us consider the fibrations
\begin{eqnarray*}
\Phi\sb{|E|} = \pr\sb 1 \circ \pi\sb P \circ \rho \sb L &:& X\to \Pl, \qquad\textrm{and} \\
\Phi\sb{|F|} = \pr\sb 2 \circ \pi\sb P \circ \rho \sb L &:& X\to \Pl,
\end{eqnarray*}
where $\pr\sb i: \Pll\to \Pl$ is the projection onto the $i$-th factor.
Because $Y\sb{(X, L)}$ has ten cusps,
the classification of fibers of (quasi-)elliptic fibrations and 
the criterion~\cite[Section 4]{RS} for quasi-ellipticity imply the following:
\begin{proposition}\label{cor:QE}
The fibrations $\Phi\sb{|E|}$ and $\Phi\sb{|F|}$
are quasi-elliptic.
Let $\Theta$ be a  fiber of the quasi-elliptic fibration $\Phi\sb{|E|}$.
Then
$\Theta$ is either of  type $\two$,
 of type $\four$ 
or of type $\fourstar$.
%
%
Moreover, we have 
\begin{eqnarray*}\label{eq:QE}
\textrm {$\Theta$ is of  type $\two$} 
&\Longleftrightarrow& \textrm{$\rho\sb L (\Theta)$ does not pass through any cusps of $Y\sb\polX$}, \\
\textrm {$\Theta$ is of  type $\four$} 
&\Longleftrightarrow& \textrm{$\rho\sb L (\Theta)$  passes through exactly one  cusp of $Y\sb\polX$}, \\
\textrm {$\Theta$ is of  type $\fourstar$} 
&\Longleftrightarrow& \textrm{$\rho\sb L (\Theta)$ is a line with multiplicity $3$ passing through} \\
&& \textrm{exactly three  cusps of $Y\sb\polX$}.
\end{eqnarray*}
Same hold for fibers of $\Phi\sb{|F|}$.
\end{proposition}
%
%
%
%
%
\begin{proof}[Proof of Theorem~\ref{thm:insep}]
Let $X$ be a supersingular $K3$ surface with $\sigma (X)\le 6$.
We choose a subcode $\CCC$ of the isotropic admissible code $\CCC\sb 7$
in Table~\ref{table:sevencodes7}
with 
$$
\dim \CCC=6-\sigma (X),
$$
and consider the corresponding overlattice $N\sb{\CCC}$ of $\Z [10\At]\oplus U(3)$.
There exists an isometry 
$$
\phi : N\sb{\CCC} \isom \NS (X)
$$
such that $\phi (e+f)$ is the class $[L]$ of a line bundle $L$ that is very ample
modulo $(-2)$-curves, where $e, f$ form the canonical basis of $U(3)$; see the proof of Theorem 1.5.
Then $Y\sb\polX$ is a complete intersection in $\P\sp 4$ with multi-degree $(2, 3)$
that has ten cusps as its only singularities.
We will prove Theorem~\ref{thm:insep} by showing that,
for this polarized supersingular $K3$ surface $(X, L)$,
the morphism $\pi\sb P$ from $Y\sb{(X, L)}$ to $\Pll$ is purely inseparable;
that is, the polynomials $a$ and $b$  in~\eqref{eq:F} are zero.
\par
\medskip
We assume that $\pi\sb P$ is separable, and derive a contradiction.
\par
\medskip
For $i=1, \dots, 10$,
let $C\sb i$ and $ D\sb i $ be the $(-2)$-curves 
contracted by $\rho\sb L$
satisfying 
$$
C\sb i ^2= D\sb i ^2=-2, 
\qquad 
C\sb i D\sb i =1,
\qquad 
\langle  [C\sb i], [D\sb i] \rangle   \perp \langle  [C\sb j], [D\sb j]\rangle \quad (i\ne j), 
$$
and let $E$, $F$ be  divisors 
such that  $\phi (e)=[E]$ and $\phi (f)=[F]$.
Then $E$ and $F$ satisfy $[E], [F]\in R\sb\polX\sperp$
and~\eqref{eq:EF}.
We put
\begin{eqnarray*}
\gamma\sb i  &:=& ([C\sb i] + 2[D\sb i])/3\;\bmod  (R\sb{\polX}\oplus R\sb{\polX}\sperp),\\
\bar f \dual &:=& [E]/ 3 \;\bmod  (R\sb{\polX}\oplus R\sb{\polX}\sperp),\\
\bar e \dual &:=& [F]/ 3 \;\bmod  (R\sb{\polX}\oplus R\sb{\polX}\sperp).
\end{eqnarray*}
The code $\CCC\sb{(X, L)}$ defined  by 
$$
\CCC\sb{(X, L)}:= \NS (X)/ (R\sb{\polX}\oplus R\sb{\polX}\sperp) \;\;\st \;\;
\Disc{R\sb{\polX}\oplus R\sb{\polX}\sperp}\cong \F\sb 3\sp{10}\oplus \F\sb 3\sp 2
$$
is isomorphic to the subcode $\CCC$ of $\CCC\sb 7$ chosen above.
Let $G$ be a divisor  on $X$.
Then $[G]\in \NS (X)$ is written as
$$
\frac{1}{3}\sum\sb{i=1}\sp{10} (s\sb i [C\sb i] + t\sb i [D\sb i]) +\frac{\alpha}{3} [E] + \frac{\beta}{3} [F],
$$
where $s\sb i, t\sb i, \alpha, \beta$ are integers satisfying $s\sb i+t\sb i\equiv 0 \bmod 3$.
We denote by
$$
\wdc{G}:= [G] \;\bmod  (R\sb{\polX}\oplus R\sb{\polX}\sperp)
$$
the word of $\CCC\sb{(X, L)}$ corresponding to $[G]$,
which is written  as
$$
(\vx( G), \bar\alpha, \bar\beta)=
\sum \sb{i=1}\sp{10} x\sb i \gamma\sb i + \bar\alpha \bar f \dual + \bar \beta \bar e \dual,
$$
where $\bar\alpha = \alpha \bmod 3$, $\bar\beta=\beta \bmod 3$, and
$$
x\sb i=\begin{cases}
0 & \textrm{if $(s\sb i, t\sb i) \equiv (0, 0)\;\bmod 3$,} \\
1 & \textrm{if $(s\sb i, t\sb i) \equiv (1, 2)\;\bmod 3$,} \\
2 & \textrm{if $(s\sb i, t\sb i) \equiv (2, 1)\;\bmod 3$.} 
\end{cases}
$$
We put
\begin{eqnarray*}
s(G) &:=& \set{i}{ (s\sb i, t\sb i)\ne (0,0)}=\set{i}{ C\sb i G \ne 0 \;\;\textrm{or}\;\; D\sb i G\ne 0}, \\ 
s\sb 1(\vx( G))&:=& \set{i}{ x\sb i \ne 0}=\set{i}{(s\sb i, t\sb i)\not\equiv  (0,0) \;\bmod 3 }, \\ 
s\sb 2(G) &:=& \set{i}{ (s\sb i, t\sb i)\ne (0,0)  \;\;\textrm{and}\;\; (s\sb i, t\sb i)\equiv  (0,0) \;\bmod 3 }.
\end{eqnarray*}
By definition, we have
$$
s (G)=s\sb 1 (\vx ( G)) \sqcup s\sb 2 (G).
$$
%
%
\begin{lemma}\label{lem:s2}
Suppose that $G$ is a reduced irreducible curve on $X$.
Then the following holds:
\begin{equation}\label{eq:Gsq}
|s\sb 2 (G)| \le \frac{1}{3} (\alpha \beta -|s\sb 1 (\vx (G))|)+1.
\end{equation}
In particular, we have $\alpha \beta -|s\sb 1 (\vx (G))| \ge -3$.
\end{lemma}
\begin{proof}
Let $s$ and $t$ be integers such that $s+t\equiv 0\;\bmod 3$.
If $(s, t)\ne (0,0)$, then
$$
\left( \nfrac{s C\sb i + t D\sb i }{3}\right)^2\le -2/3
$$
holds. If $(s, t)\ne (0,0)$ and $(s, t)\equiv (0,0)\;\bmod 3$, then
$$
\left( \nfrac{s C\sb i + t D\sb i }{3}\right)^2\le -2
$$
holds. Therefore we have
$$
G^2\le -\frac{2}{3} |s\sb 1 (\vx (G))| - 2 |s\sb 2 (G)| +\frac{2}{3} \alpha\beta.
$$
On the other hand, we have $G^2\ge -2$.
Hence we get the inequality~\eqref{eq:Gsq}.
\end{proof}
Let us denote by $\bar T$ the Cartier divisor on $Y\sb{\polX}$ cut out by the equation
\begin{equation}\label{eq:T}
\frac{\partial \Psi}{\partial W}=-aW+b=0,
\end{equation}
and let $T$ be the proper transform of $\bar T$ by $\rho \sb L$.
By the assumption that $\pi\sb P$ is separable,
$\bar T$ is a divisor and 
 $\pi\sb P$ is \'etale outside $\bar T$.
Hence the divisor $\bar T$ contains the ten cusps of $Y\sb{(X, L)}$.
Therefore  we have 
\begin{equation}\label{eq:sT}
s(T)=\{ 1, 2, \dots, 10\}.
\end{equation}
From the defining equation~\eqref{eq:T} of $\bar T$ on $Y\sb{(X, L)}$, we have
\begin{equation}\label{eq:ETFT}
ET=FT=6.
\end{equation}
We denote by $C\sb E$  the closure of the locus 
$$
\set{x\in X}{\textrm{the fiber of $\Phi\sb{|E|}$ passing through $x$ is of type $\two$ and is singular at $x$}},
$$
and equip $C\sb E$ with the reduced structure.
A general member $E$ of $|E|$ intersects $C\sb E$ at one point with multiplicity $3$~(\cite{BM}).
See Figure~\ref{fig:CE}.
%
%
%
%
We define $C\sb F$ in the same way.
Both of $C\sb E$ and $C\sb F$ are irreducible,
and 
 we have
\begin{equation}\label{eq:CEECFF}
C\sb E E =C\sb F F=3.
\end{equation}
Because $\pr\sb 1 :\Pll\to\Pl$ is smooth,
if $\pr\sb 1\circ \pi\sb P$ is not smooth at a non-singular point of $Y\sb\polX$,
then $\pi\sb P$ is not smooth at that point.
Therefore
the divisor  $T$ contains $C\sb E$ as a reduced  irreducible component.
Same holds for $C\sb F$.
\begin{claim}\label{claim:distinct}
The two curves $C\sb E$ and $C\sb F$ are distinct.
\end{claim}
\begin{proof}
Suppose that $C\sb E=C\sb F$ holds.
Let $x$ be a general point of $C\sb E=C\sb F$.
Since the fibers $E\sb x$ of $\Phi\sb{|E|}$ and $F\sb x$ of $\Phi\sb{|F|}$ 
passing through $x$ are both singular at $x$,
we have $E\sb x F\sb x\ge 4$,
which contradicts $EF=3$.
\end{proof}
Let 
$$
T=C\sb E + C\sb F + T\sb 1+\cdots+T\sb{t}
$$
be the decomposition of $T$ into reduced irreducible components.
We put
\begin{eqnarray*}
[C\sb E] &=& \sum (s\sb{E,i} [C\sb i] + t\sb {E, i} [D\sb i] )/3 + (\alpha \sb E [E] +\beta\sb E [F])/3,\cr
[C\sb F] &=& \sum (s\sb{F,i} [C\sb i] + t\sb {F, i} [D\sb i] )/3 + (\alpha \sb F [E] +\beta\sb F [F])/3, \cr
[T\sb \nu] &=&  \sum (s\sb{\nu,i} [C\sb i] + t\sb {\nu , i} [D\sb i] )/3 + (\alpha \sb \nu [E] +\beta\sb \nu [F])/3
\quad (\nu=1, \dots, t) .
\end{eqnarray*}
Since $E$ and $F$ are nef,
we have
\begin{equation}\label{eq:nef}
\alpha \sb E\ge 0,\;\;
\beta \sb E\ge 0,\;\;
\alpha \sb F\ge 0,\;\;
\beta \sb F\ge 0,\;\;
\alpha \sb \nu\ge 0,\;\;
\beta \sb \nu\ge 0\quad (\nu=1, \dots, t).
\end{equation}
Since $\pi\sb P$ is  finite, $\pi\sb P\circ \rho\sb L$ maps each irreducible component of $T$ to a curve on $\Pll$.
Therefore we have 
\begin{equation}\label{eq:finite}
\alpha \sb \nu> 0\;\;\textrm{or}\;\;
\beta \sb \nu> 0\qquad (\nu=1, \dots, t).
\end{equation}
By~\eqref{eq:CEECFF}, we have
\begin{equation}\label{eq:33}
\beta\sb E=3, \quad \alpha\sb F=3.
\end{equation}
Then, from~\eqref{eq:ETFT}, we have 
\begin{equation}\label{eq:sum3}
\alpha\sb E +\sum\sb{\nu=1}\sp t \alpha\sb{\nu}=3, \quad\textrm{and}\quad \beta\sb F+\sum\sb{\nu=1}\sp t \beta\sb{\nu}=3.
\end{equation}
Consider the words
$$
\wdc{C\sb E}=(\vx\sb E, \bar{\alpha}\sb E, 0), \quad
\wdc{C\sb F}=(\vx\sb F, 0, \bar{\beta}\sb F), \quad
\wdc{T\sb{\nu}}=(\vx\sb {\nu}, \bar{\alpha}\sb \nu , \bar{\beta}\sb \nu)\qquad (\nu=1, \dots, t)
$$
in the code $\CCC\sb{(X, L)}$.
From Lemma~\ref{lem:s2}, we have 
\begin{equation}\label{eq:effv}
-\wt(\vx\sb\nu) + \alpha \sb \nu  \beta\sb \nu  \ge -3\qquad (\nu=1, \dots, t).
\end{equation}
\begin{claim}\label{eq:EFzero}
$\vx\sb E=\vx\sb F=\vz$.
\end{claim}
\begin{proof}
Let $\Theta$ be a  fiber of $\Phi\sb{|E|}$ such that
$\rho\sb L (\Theta)$ passes  through a  cusp $q\sb i:=\rho\sb L (C\sb i)=\rho\sb L (D\sb i)$
of $Y\sb{\polX}$.
Then $\Theta$ is of type $\four$ or $\fourstar$.
Suppose that $\Theta$ is of type $\four$.
Then $\Theta$ consists of three irreducible components of multiplicity one,
two of which are $C\sb i$ and $D\sb i$,
that intersect at one point.
The curve $C\sb E$ passes through the intersection point.
Since $\Theta C\sb E=3$,
we have $C\sb E C\sb i = C\sb E D\sb i =1$,
and therefore $s\sb {i, E}= t\sb{i, E}=-3$ holds.
Suppose that $\Theta$ is of type $\fourstar$.
Then $C\sb E$ passes through a point of the multiplicity $3$ component  
of $\Theta$, and does not intersect other irreducible components.
This fact can be proved by considering the pull-back of the quasi-elliptic fibration $\Phi\sb{|E|}$
by the base change $\Pl\to\Pl$ of degree $2$ branching at the point  $\Phi\sb{|E|} (\Theta)$,
which makes the fiber $\Theta$ into type $\four$.
Then it follows that $s\sb {i, E}=t\sb {i, E}=0$.
In any case,
we have  $i\notin s\sb 1 (\vx\sb E)$.
Since this holds for any cusp $q\sb i$ of $Y\sb{\polX}$,
we have $s\sb 1 (\vx\sb E)=\emptyset$.
\end{proof}
Since $\bar T$ is a Cartier divisor of $Y\sb\polX$,
the total transform  $\rho\sb L\sp *({\bar T}) $  is contained in $R\sb{\polX}\sperp=\Z[E]\oplus \Z [F]$,
and hence   $\wdc{T}=0$. Therefore we obtain 
\begin{equation}\label{eq:Cartier}
 \vx\sb 1 +\cdots +\vx\sb {t} =\vz.
\end{equation}
By~\eqref{eq:sT},
we have
$$
 s(C\sb E) \cup s (C\sb F) \cup s (T\sb{1}) \cup \dots \cup s (T\sb {t})=\{1, 2, \dots, 10\}.
$$
Using Lemma~\ref{lem:s2},
we obtain
\begin{multline}\label{eq:theineq}
t+2+
\frac{1}{3} \Bigl (\alpha\sb E\beta\sb E + \alpha\sb F \beta\sb F +
 \sum\sb{\nu=1}\sp{t} (\alpha\sb \nu \beta\sb \nu - |s\sb 1 (\vx \sb  \nu )|)\Bigr) \ge\\ 
10 - |s\sb 1 (\vx  \sb{1}) \cup \dots \cup 
s\sb 1 (\vx \sb {t})|.
\end{multline}
Because $\CCC\sb{(X, L)}$ is isomorphic to a subcode of $\CCC\sb 7$, 
 we have shown that 
 there exist  integers $\alpha\sb E$, $\beta\sb F$, $\alpha\sb \nu, \beta\sb \nu$ $(\nu=1, \dots, t)$
and 
words
$$
(\vz, \bar{\alpha}\sb E, 0), \quad
(\vz, 0, \bar{\beta}\sb F), \quad
(\vx\sb {\nu}, \bar{\alpha}\sb \nu , \bar{\beta}\sb \nu)\qquad (\nu=1, \dots, t)
$$
in the code $\CCC\sb 7$ satisfying~\eqref{eq:nef}-\eqref{eq:theineq}.
Using computer, however, we can show that such integers and words do not exist.
Thus we get a contradiction.
\end{proof}
Instead of the code $\CCC\sb 7$,
we can use the codes $\CCC\sb 3, \dots, \CCC\sb 6$ in Table~\ref{table:sevencodes7}.
However, we cannot use $\CCC\sb 2$ or $\CCC\sb 1$.
Indeed, in $\CCC\sb 1$,
for example,
we have the following integers and words:
\begin{eqnarray*}
&&[0, 0, 0, 0, 0, 0, 0, 0, 0, 0, 3, 0], \\
&&[0, 0, 0, 0, 0, 0, 0, 0, 0, 0, 0, 3], \\
&&[1, 0, 0, 0, 0, 0, 0, 0, 1, 1, 0, 1], \\
&&[0, 2, 0, 0, 0, 0, 0, 2, 0, 2, 1, 0], \\
&&[0, 0, 2, 0, 0, 0, 2, 0, 2, 0, 1, 0], \\
&&[0, 0, 0, 1, 0, 0, 1, 1, 0, 0, 0, 1], \\
&&[2, 1, 1, 2, 0, 0, 0, 0, 0, 0, 1, 1].
\end{eqnarray*}
Nevertheless, we can ask the following:
\begin{question}\label{question1}
Is $\pi\sb P: Y\sb{\polX}\to\Pll$ inseparable
for any polarized supersingular $K3$ surface $(X, L)$ of degree $6$ with $\RRR \sb{(X, L)}= 10 A\sb 2$
and $R\sperp\sb{(X, L)}\cong U(3)$?
\end{question}
\begin{example}\label{example:sigmaG0}
Consider the purely inseparable triple cover of $\Pll$ defined by
$$
W^3=(x^3-x)(y^3-y),
$$
and the corresponding polarized supersingular $K3$ surface $(X, L)$.
We will show that the Artin invariant of $X$ is $1$, and that the $5$-dimensional ternary code
$\CCC\sb{\polX}$ is isomorphic to $\CCC\sb 1$.
For $\alpha\in \F\sb 3$,
let $l\sb\alpha$ and $m\sb\alpha$ be the lines on $\Pll$
defined by 
$x=\alpha$ and $y=\alpha$, respectively.
The strict transforms of $l\sb\alpha$ and $m\sb\alpha$ by $\pi\sb P\circ \rho\sb L$
are written as $3 \tilde l\sb{\alpha}$ and $3\tilde m\sb\alpha$,
respectively.
Numbering  the twenty $(-2)$-curves $C\sb 1, D\sb 1, \dots, C\sb{10}, D\sb{10}$
in an appropriate way,
we can write the numerical equivalence classes 
$[\tilde l\sb{\alpha}]$,
$[\tilde m\sb{\alpha}]$
as follows:
\begin{eqnarray*}
[\tilde l\sb{0}] &=& A\sb 1 + A\sb 2 + A\sb 3 + [E]/3, \cr
[\tilde m\sb{0}] &=& A\sprime\sb 1 + A\sprime\sb 4 + A\sprime\sb 7 + [F]/3, \cr
[\tilde l\sb{1}] &=& A\sb 4 + A\sb 5 + A\sb 6 + [E]/3, \cr
[\tilde m\sb{1}] &=& A\sprime\sb 2 + A\sprime\sb 5 + A\sprime\sb 8 + [F]/3, \cr
[\tilde l\sb{2}] &=& A\sb 7 + A\sb 8 + A\sb 9 + [E]/3, \cr
[\tilde m\sb{2}] &=& A\sprime\sb 3 + A\sprime\sb 6 + A\sprime\sb 9 + [F]/3, 
\end{eqnarray*}
where
$$
A\sb i= -([C\sb i] + 2 [D \sb i])/3,
\qquad 
A\sprime \sb i = -(2 [C\sb i] + [D \sb i])/3.
$$
The discriminant of the sublattice of $\NS (X)$ generated by the classes  $[E], [F]$, the classes of  
the twenty exceptional curves, and 
the $6$ classes above 
is equal to $-9$.
Hence these classes span $\NS (X)$,
and the Artin invariant of $X$ is $1$.
The $6$ words 
$\wdc{\tilde l\sb{\alpha}}$,
$\wdc{\tilde m\sb{\alpha}}$
generate a $5$-dimensional ternary code isomorphic to $\CCC\sb 1$.
\end{example}
\begin{question}\label{question2}
Find the defining equations 
of  purely inseparable triple covers of $Q=\Pll$ 
corresponding to the other ternary codes $\CCC\sb 2, \dots, \CCC\sb 7$
of dimension $5$
 in Table~\ref{table:sevencodes7}.
(See Corollary~\ref{cor:seven}.)
\end{question}
In~\cite{DK},
Dolgachev and Kondo gave various defining equations of the supersingular $K3$ surface
in characteristic $2$ with Artin invariant $1$,
and determined the full automorphism group of this $K3$ surface.
We expect that  various defining equations of 
the supersingular $K3$ surface
in characteristic $3$ with Artin invariant $1$
would be also helpful  in the study of  the automorphism group of this surface.


\bibliographystyle{amsplain}

\providecommand{\bysame}{\leavevmode\hbox to3em{\hrulefill}\thinspace}

\end{document}